\documentclass[11pt]{article}
\usepackage{graphicx}
\usepackage{amsmath}
\usepackage{latexsym, color}
\def\bega{\begin{array}}
\def\enda{\end{array}}

\def\begi{\begin{itemize}}
\def\endi{\end{itemize}}

\def\forall{\hbox{for all}~}
\def\L{{\bf L}}
\def\ds{\displaystyle}

\def\argmax{\hbox{arg}\!\max}

\def\wto{\rightharpoonup}
\def\argmin{\hbox{arg}\!\min}
\def\ve{\varepsilon}

\def\I{{\cal I}}

\def\R{I\!\!R}

\def\Supp{\hbox{Supp}}
\def\implies{\Longrightarrow}
\def\vp{\varphi}

\def\G{{\cal G}}

\def\v{\vskip 1em}
\def\O{{\cal O}}
\def\C{{\cal C}}

\def\ov{\overline}

\def\Hat{\widehat}

\def\bel{\begin{equation}\label}
\def\eeq{\end{equation}}
\def\sqr#1#2{\vbox{\hrule height .#2pt
\hbox{\vrule width .#2pt height #1pt \kern #1pt
\vrule width .#2pt}\hrule height .#2pt }}
\def\square{\sqr74}
\def\endproof{\hphantom{MM}\hfill\llap{$\square$}\goodbreak}

\newtheorem{theorem}{Theorem}[section]
\newtheorem{lemma}[theorem]{Lemma}

\newtheorem{definition}[theorem]{Definition}
\newtheorem{remark}[theorem]{Remark}

\begin{document}
\title{\bf Globally Optimal Departure Rates  for  Several Groups of Drivers}

\author{Alberto Bressan$^*$ and Yucong Huang$^{**}$\\
\,\\
\small (*) Department of Mathematics, Penn State University,
University Park, PA 16802, U.S.A.
\,\\ 
\small  (**)  Mathematical Institute,
University of Oxford, Woodstock Road, Oxford OX2 6GG, UK. 
\\ \, 
\small
E-mails:  axb62@psu.edu, ~yucong.huang@maths.ox.ac.uk }
\maketitle

\begin{abstract} The first part of this paper contains a brief introduction to 
conservation law models of traffic flow on a network of roads. 
Globally optimal solutions and Nash equilibrium solutions are reviewed,
with several groups of drivers sharing different cost functions.

In the second part we consider  a globally optimal set of departure rates, for
different groups of drivers but on a single road.
Necessary conditions are proved, 
which lead to a practical algorithm for computing the optimal solution.
\end{abstract}

\section{Introduction}
\label{sec:0}
\setcounter{equation}{0}

Macroscopic models of traffic flow, first introduced in 
\cite{LW, R}, have now become a topic of extensive research. 
On a single road, the evolution of the traffic density can be described by
a scalar conservation law.  In order to extend the model to a whole network of
roads, additional boundary conditions must be inserted, 
describing traffic flow at each intersection; see \cite{BN1, CGP, GHP, GP, HMR, HR}
or the survey \cite{BCGHP}.
A major eventual goal of these models is to understand traffic patterns,
determined by the behavior of a large number of drivers
with different origins and destinations.

In a basic  setting,  one can consider $N$ groups of drivers, say $\G_1,\ldots, \G_N$.
Drivers from each group have the same origin and destination, and a cost
which depends on their departure and arrival time.  For such a model, 
two kind of solutions are of interest:
\begi
\item[-] The Nash equilibrium solution, where each driver chooses his own departure time
and route to destination, in order to minimize his own cost.

\item[-] The global optimization problem, where a central planner seeks to schedule
all departures in order to minimize the sum of all costs.
\endi
In general, these criteria determine very different traffic patterns. 
To fix the ideas, let $t\mapsto u_i(t)$, $i=1,\ldots, N$,
be the departure rate of drivers of the $i$-th group, so that
$$\int_{-\infty}^t u_i(s)\, ds$$
yields the total number of these drivers who depart before time $t$.
We recall that the support of $u_i$, denoted by
Supp$(u_i)$, is  the closure of set of times $t$ where $u_i(t)>0$.

Roughly speaking, the two above solutions can be characterized
as follows.

{\bf (I)}  In a Nash equilibrium, all drivers within the same group pay the same cost.
Namely, there exists constants $K_1,\ldots, K_N$ such that 
\begi
\item every driver of the $i$-th group, departing at a time $t\in Supp(u_i)$ 
bears the cost $K_i$.

\item if a driver of the $i$-th group were to depart at any time $t\in \R$
(possibly outside the support of $u_i$), 
he would incur in a cost $\geq K_i$. 
\endi

{\bf (II)}  For a global optima, there exist constants $C_1,\ldots, C_n$
(where $C_i$ is the marginal cost for adding one more driver of the $i$-th group) 
such that 
\begi
\item If one additional driver of the $i$-th group is added at any time $t\in \Supp(u_i)$,
then the total cost increases by $C_i$.

\item  If one additional driver of the $i$-th group is added at any time $t\in \R$
(possibly outside the support of $u_i$),
then the increase in the total cost is greater or equal to $C_i$.
\endi

\begin{figure}[htbp]
\centering
  \includegraphics[width=0.9\textwidth]{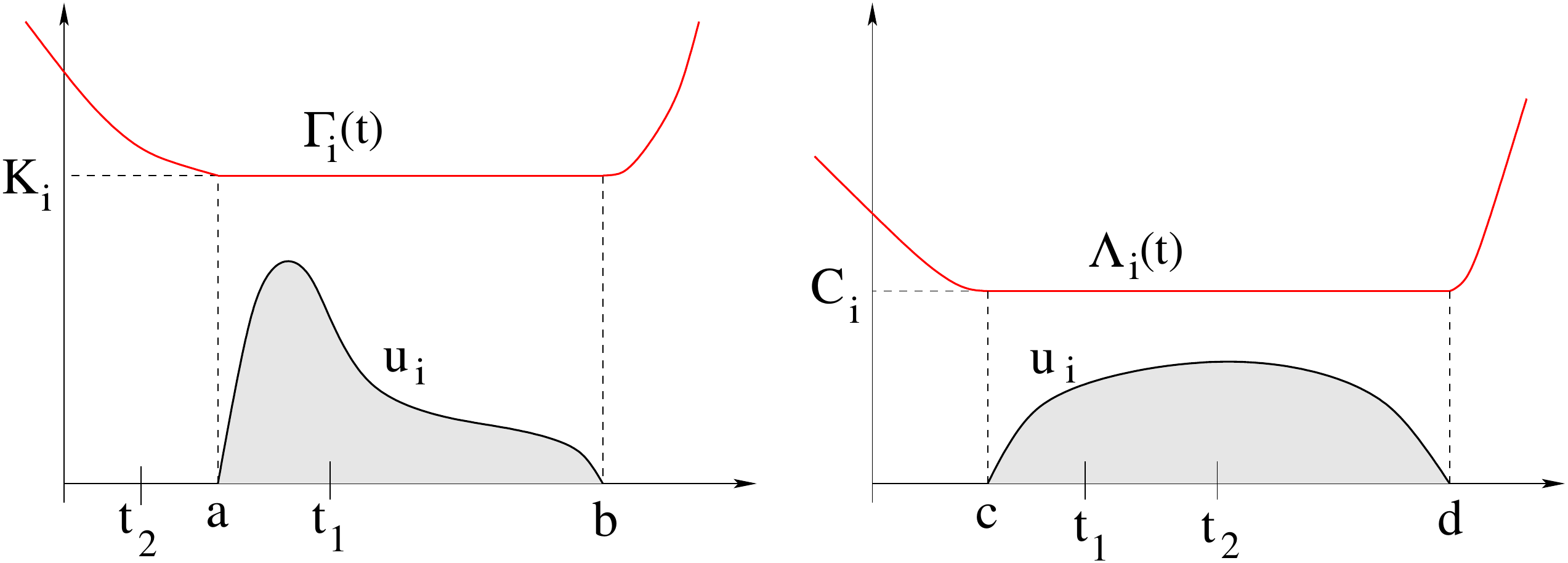}
    \caption{\small  Left:  the rate of departures $u_i$, for drivers of the $i$-th group,
    in a Nash equilibrium solution.
Each of these drivers starts at some time $t\in [a,b]$. To achieve an equilibrium,
    the cost $\Gamma_i(t)$ to any driver departing at time $t$ must be constant inside 
    $[a,b]$ and larger outside.  
    Right:  the rate of departures $u_i$, in a globally optimal solution.
   Here $\Lambda_i(t)$ denotes the marginal  cost for
   inserting an additional driver of the $i$-th group, departing at time $t$.
To achieve global optimality, $\Lambda_i$ must be constant on the support of $u_i$, and larger outside.}
\label{f:tf238}
\end{figure}

At an intuitive level, these conditions are easy to explain.
In Fig.~\ref{f:tf238}, left, the function $\Gamma_i(t)$ 
denotes the cost to an $i$-driver departing at time $t$.
  If $\Gamma_i$ did not attain
its global minimum simultaneously at all points $t\in [a,b]= \Supp(u_i)$,
 then we could find times $t_1\in [a,b]$
and $t_2\in\R$ such that 
$\Gamma_i(t_2)<\Gamma_i(t_1)$.  In this case,  the driver departing at 
time $t_1$ could lower his own cost choosing to depart at time $t_2$ instead.
This contradicts the definition of equilibrium.

In Fig.~\ref{f:tf238}, right, the function $\Lambda_i(t)$ denotes the marginal cost
for inserting one additional driver of the $i$-th family, departing at time $t$.  
This accounts for the additional cost to the new driver, and also for the increase
in the cost to all other drivers who are slowed down by the presence of
one more car on the road.
If $\Lambda_i$ did not attain its global minimum at all points in $[c,d]= \Supp(u_i)$,
then we could find times $t_1\in [c,d]$
and $t_2\in\R$ such that 
$\Lambda_i(t_2)<\Lambda_i(t_1)$.  
In this case we could consider a new traffic pattern, with
 one less driver departing at time $t_1$ and one more departing at time $t_2$.
 This would achieve a smaller total cost, contradicting  the 
 assumption of optimality.

While the criterion {\bf (I)} for an equilibrium solution is easy to justify,
a rigorous proof of the necessary condition {\bf (II)}  for a global optimum
faces considerable difficulties. Indeed, to compute the ``marginal cost"
for adding one more driver, one should differentiate the solution of a
conservation law w.r.t.~the initial data (or the boundary data).   As it is well known, 
in general one does not have enough regularity to carry out such a differentiation.
To cope with this difficulty
one can introduce  a ``shift differential", describing how the shock
locations change, depending on parameters. See \cite{BM1, BM2, BS, PU, U}
for results in this direction.

The first part of this paper contains an introduction to macroscopic models of
traffic flow on a network of roads.  Section~\ref{s:2}  starts by reviewing the classical
LWR model for traffic flow on a single road, in terms of a scalar conservation law for
the traffic density.   We then discuss various boundary conditions, modeling traffic flow
at an intersection.  
Finally, given a cost function depending on the departure and arrival times of each driver,
we review the concepts of globally optimal solution and of  Nash equilibrium solution.

The second part of paper contains original results.   We consider here 
$N$ groups of drivers traveling along the same road, but with different
departure and arrival costs.   We seek departure rates $u_1(\cdot),\ldots, u_N(\cdot)$
which are globally optimal. Namely, they minimize the sum of all costs to all drivers.
A set of necessary conditions for optimality is derived, thus extending
the result in \cite{BH1} to the case where several groups of drivers are present.
Relying on these conditions, in the last section we introduce an algorithm that 
numerically computes such globally optimal solutions.

For an introduction to the general theory of conservation laws we refer to
\cite{Bbook, Evans, Smoller}.
A more comprehensive discussion of  various models of traffic flow
can be found in  \cite{Bellomo, BD, Dag, GHP}.

\section{Conservation law models for traffic flow}
\label{s:2}
\setcounter{equation}{0}
\subsection{Traffic flow on a single road.}
According to the classical LWR model \cite{LW, R}, traffic density on 
a single road can be described in terms of a scalar conservation law
\bel{claw}
\rho_t(t,x) + f(\rho(t,x))_x~=~0.\eeq
Here
$t$ is the time, while $x\in\R$ is the space variable along the road. Moreover
 \begi
\item $\rho$ is the {\bf traffic density}, i.e., the number of cars per unit length of the road.

\item $v= v(\rho)$ is the {\bf velocity of cars}, which we assume depends only  the traffic density.

\begin{figure}[htbp]
\centering
  \includegraphics[width=0.9\textwidth]{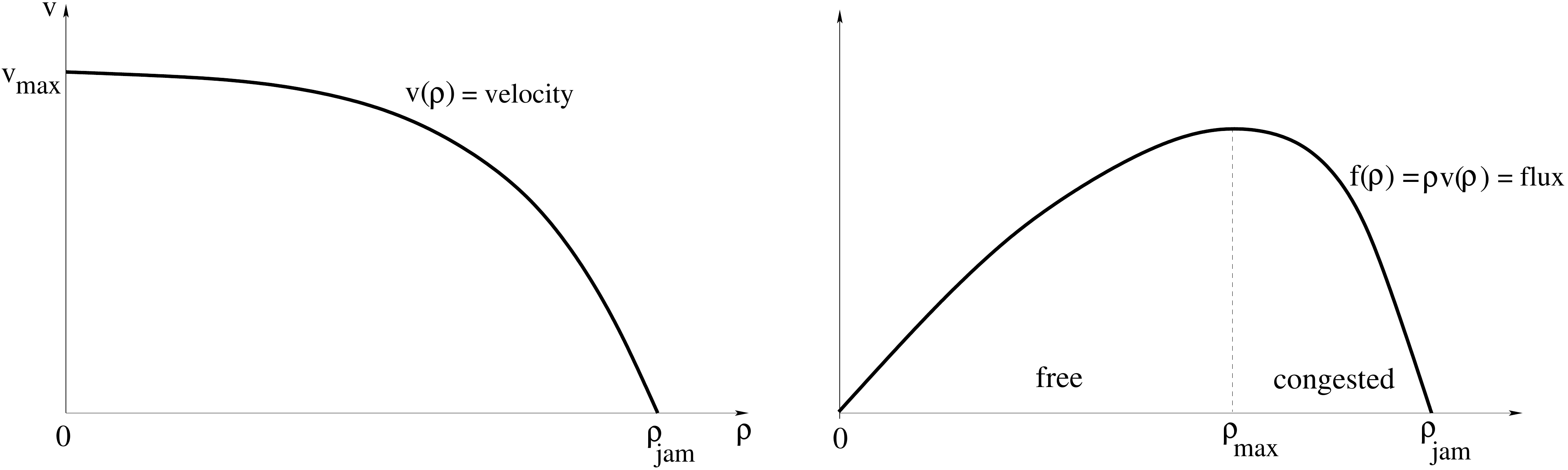}
    \caption{\small  Left: the velocity of cars as a function of the traffic density.
 Notice that $v$ is maximum when $\rho=0$ and the road is empty. 
 The velocity decreases to zero  as the density approaches  a critical density 
 $\rho_{jam}$, where cars are packed bumper-to-bumper and no one moves.
  Right: the flux function $f$, depending on the density. 
  Typically, this function is concave down, 
  vanishes at $\rho=0$ and at $\rho = \rho_{jam}$, and has 
  a maximum at some intermediate point $\rho_{max}$. 
  }
\label{f:tf239}
\end{figure}

\item $f= f(\rho)$ is the
 {\bf flux}, i.e., the number of cars crossing a
point $x$ along the road, per unit time. 
We have the identity
$$\hbox{[flux]}~=~\hbox{[density]$\times$[velocity]}~=~\rho\cdot v(\rho)$$
\endi
As shown in Fig.~\ref{f:tf239}, the velocity should be a decreasing 
function of the car density.  Concerning the flux function, a 
natural set of assumptions is
\bel{fi}
f\in \C^2, \qquad f''<0,\qquad f(0)= f(\rho_{jam}) = 0.\eeq
Here $\rho_{jam}$ is the maximum density of cars allowed 
on the $k$-th road.  This corresponds to bumper-to-bumper packing,
where no car can move.

Smooth solutions of the conservation law (\ref{claw})  can be computed by the classical
method of characteristics. By the chain rule, one obtains
\bel{qe} \rho_t + f'(\rho) \rho_x~=~0.\eeq
Hence, if $t\mapsto x(t)$ is a curve such that
\bel{char}\dot x(t)~\doteq~{d\over dt} x(t)~=~f'(\rho(t,x(t))),\eeq
then
the equation (\ref{qe}) yields
$${d\over dt} \rho(t, x(t))~=~\rho_t+\rho_x\, \dot x ~=~0.$$
In other words, the density is constant along each characteristic curve
satisfying (\ref{char}). 
Notice that the assumptions (\ref{fi}) imply the inequality 
$$\hbox{[car speed]}~\dot=~v(\rho)~=~f(\rho)/\rho~\geq~
f'(\rho)~=~v(\rho) + \rho\, v'(\rho)~=~\hbox{[characteristic speed]}.
$$
With reference to Fig.~\ref{f:tf239},
let $\rho_{max}$ be the density at which the flux is maximum.
We say that a state $\rho$ is
\begi
\item {\bf free}, if $\rho<\rho_{max}$, hence the characteristic speed $f'(\rho)$ is positive,
\item {\bf congested},   if $\rho>\rho_{max}$, hence the characteristic speed $f'(\rho)$ is negative.
\endi
\v
Due to the non-linearity of the flux function $f$, It is well known that solutions can develop
shocks in finite time. The conservation law (\ref{claw})  must 
thus be interpreted in distributional sense. 
For the general theory
of entropy weak solutions to conservation laws, we refer to \cite{Bbook, Smoller}.

\subsection{Traffic flow at  road intersections.}

To model vehicular traffic on an entire network of roads, the conservation laws
describing traffic flow on each road must be supplemented with
boundary conditions, describing the behavior at road intersections.

\begin{figure}[htbp]
\centering
  \includegraphics[width=0.4\textwidth]{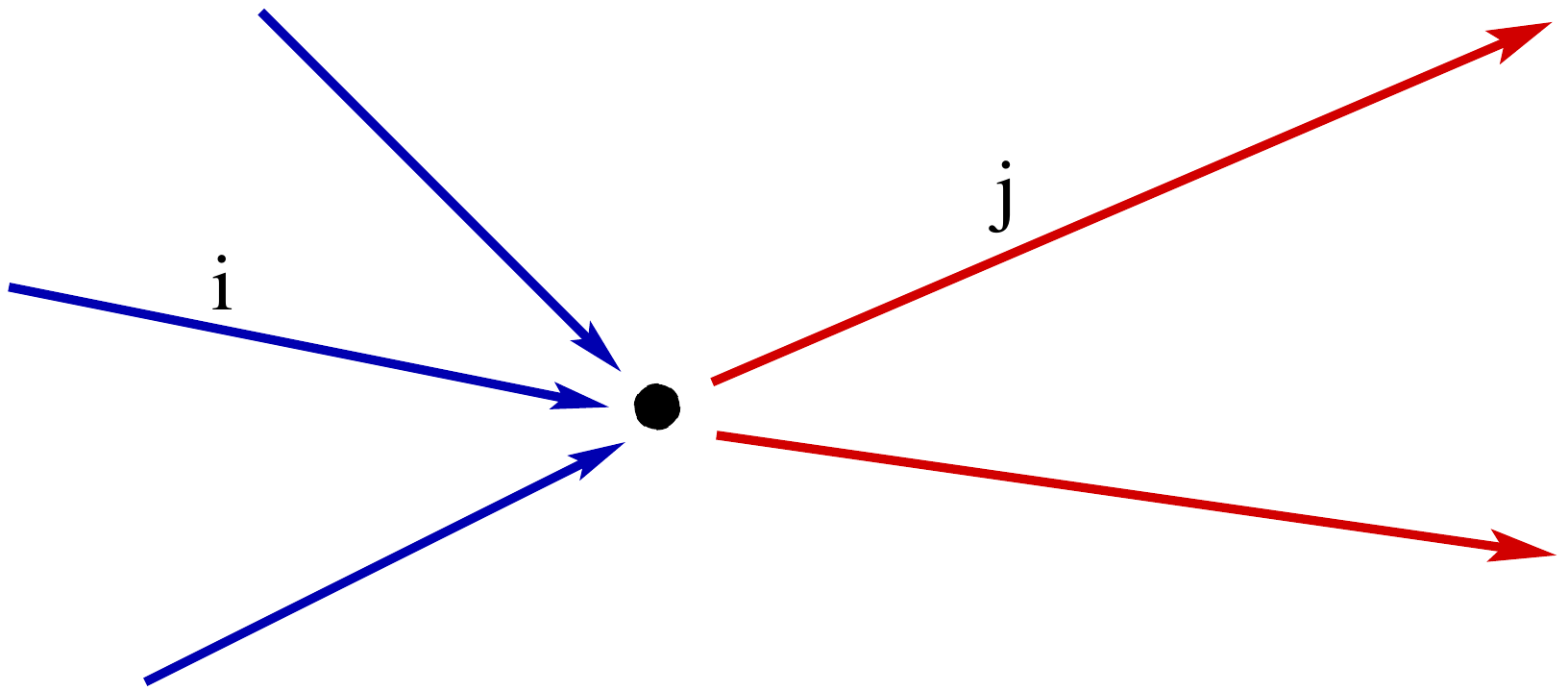}
    \caption{\small An intersection with 3 incoming and 2 outgoing roads.}
\label{f:tf41}
\end{figure}

Consider an intersection, say with $m$ incoming roads 
$i\in \{1,\ldots, m\}= \I$
and $n$ outgoing roads $j\in \{m+1,\ldots,m+n\}=\O$, see Fig.~\ref{f:tf41}.
We shall use the space variable 
$x\in \,]-\infty,\, 0]$  for incoming roads and 
$x\in [0, +\infty[\,$ for outgoing roads.
Throughout the following we assume that the density of traffic on 
each road is governed by a conservation law
\bel{clk} \rho_t + f_k(\rho)_x~=~0,\qquad\qquad f_k(\rho) \,= \,\rho \,v_k(\rho),\eeq
where the flux function $f_k$ satisfies (\ref{fi}), for every $k=1,\ldots, m+n$. 

An appropriate model must depend on various parameters, namely 
\begi
\item $c_i$ = relative priority of drivers arriving from road $i$.
\item $\theta_{ij}$  = fraction of drivers from road $i$ that turn into road $j$. 
\endi
For example, if the intersection is regulated by a crosslight, $c_i$ could
measure the fraction of time when drivers from road $i$ get green light, on average.
It is natural to assume
\bel{tca}c_i,~\theta_{ij}~\geq~0,\qquad
\qquad\sum_{i \in \I} c_i ~=~1,
\qquad\qquad  \sum_{j\in \O}\theta_{ij} ~=~1.\eeq
Boundary conditions should determine the limit
values of the traffic density on each of the $m+n$
roads meeting at the intersection:
\bel{bval} \rho_i(t, 0-)\,=\,\lim_{x\to 0-} \rho_i(t,x), \qquad \rho_j(t, 0+)\,=\,\lim_{x\to 0+} \rho_j(t,x),
\qquad\qquad\forall i\in\I,~ j\in\O.\eeq
At first sight, one might guess that $m+n$ conditions will be required.
However, this is not so, because on some roads the characteristics 
move toward the intersection.  For these roads, the limits in (\ref{bval}) are already
determined by integrating along characteristics.
Boundary conditions are required only for those roads where the characteristics 
move away from the origin.  Recalling the definition of free and congested states, 
we thus have
$$\bega{l} \hbox{[\# of boundary conditions needed to determine the flux at the intersection]}\\[3mm]
\quad ~=~\hbox{
[\# of incoming roads which are congested ] + [\# of outgoing roads which are free].}
\enda $$
It now becomes apparent that, to assign a meaningful set of  boundary conditions,  
 several different cases must be considered.   
 
 To circumvent  these difficulties,
 an alternative approach  developed by Coclite, Garavello,  
 and Piccoli \cite{CGP, GP, GP2}
relies on the construction of a {\bf Riemann Solver}.  
Instead of assigning a variable number of  boundary conditions, 
here the idea is to introduce a rule for solving all {\bf Riemann problems} 
(i.e.~the initial-value problems where at time $t=0$ the 
densities $\rho_k$ and turning preferences $\theta_{ij}$ 
are constant along each road).  
Relying on front-tracking approximations, under suitable conditions
one can prove that the solutions with general initial data 
are also uniquely determined.
\v
We briefly review the main steps of this construction, for the constant initial data
$$\left\{ \bega{l}
\rho_1,\ldots, \rho_m, \rho_{m+1},\ldots \rho_{m+n}~=~\hbox{initial densities on the incoming and outgoing roads,}\cr
\theta_{ij} ~ =~ \hbox{fraction of drivers from road $i$ that turn into road $j$.}
\enda\right.$$
\v
{\bf Step 1.}   Determine the maximum flux $f_i^{max}$ that can exit from each incoming road
$i\in\I$.

\begin{figure}[htbp]
\centering
\includegraphics[width=0.8\textwidth]{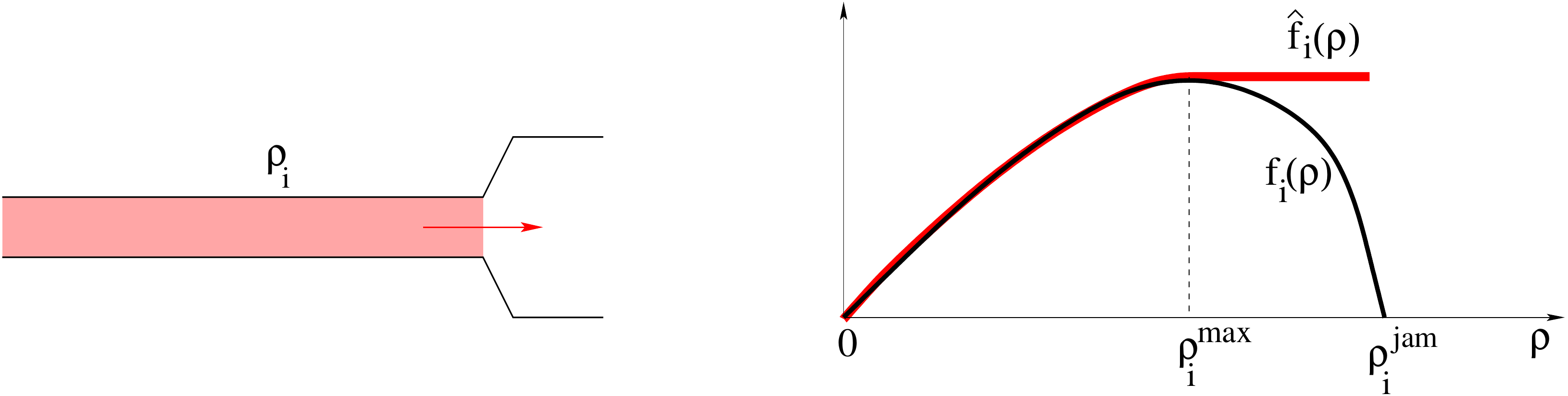}
    \caption{\small Computing the maximum flux that can come out from a road $i\in\I$.}
\label{f:tf45}
\end{figure}

As shown in Fig.~\ref{f:tf45}, this is computed by
$$f_i^{max}~=~ \Hat f_i(\rho_i)~=~\left\{ \bega{cl} f_i(\rho_i)\qquad &\hbox{if}\qquad \rho_i\leq\rho_i^{max},\\[3mm]
 f_i(\rho_i^{max})\qquad &\hbox{if}\qquad \rho_i>\rho_i^{max}\,.\enda\right.$$

%$$f_i^{max}~=~\left\{ \bega{cl} f_i(\rho_i)\quad& \hbox{if}\quad \rho_i\leq \rho_i^{max}~~~\hbox{(free state)}\\[3mm]
% f_i(\rho_i^{max})\quad& \hbox{if}\quad \rho_i\geq \rho_i^{max}~~~\hbox{(congested state)}\enda\right.$$
\v
{\bf Step 2:}   Determine the maximum flux $f_j^{max}$ that can enter each outgoing road
$j\in\O$.

\begin{figure}[htbp]
\centering
\includegraphics[width=0.8\textwidth]{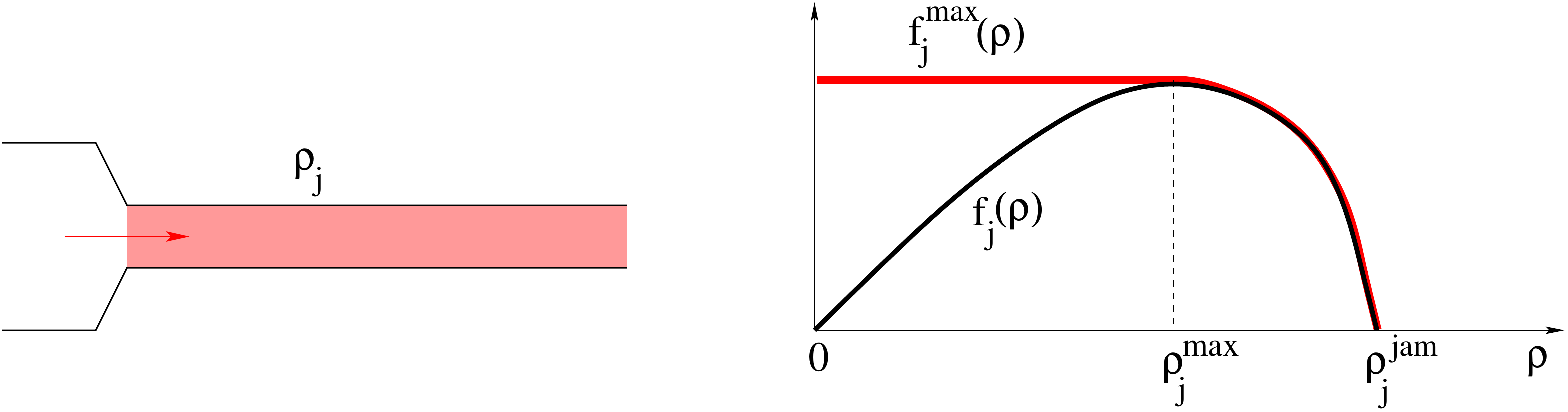}
    \caption{\small Computing the maximum flux that can  get into a road $j\in\O$.}
\label{f:tf46}
\end{figure}

As shown in Fig.~\ref{f:tf46}, this is computed by
$$f_j^{max}~=~ \Hat f_j(\rho_j)~=~\left\{ \bega{cl} f_j(\rho_j)\qquad &\hbox{if}\qquad \rho_i\geq\rho_i^{max},\\[3mm]
 f_j(\rho_j^{max})\qquad &\hbox{if}\qquad \rho_j<\rho_j^{max}\,.\enda\right.$$

{\bf Step 3:}  Given the maximum incoming and outgoing fluxes $f_i^{max}$,    $f_j^{max}$, and the turning preferences $\theta_{ij}$, determine the  region
of {\bf admissible incoming fluxes}
\bel{Omdef}
\Omega~\doteq~\left\{ (f_1,\ldots, f_m)\,;~~f_i\in [0, f_i^{max}], \qquad
\sum_{i\in\I} f_i \theta_{ij}\leq f_j^{max}\quad\forall j\in\O\right\}.\eeq

\begin{figure}[htbp]
\centering
\includegraphics[width=0.7\textwidth]{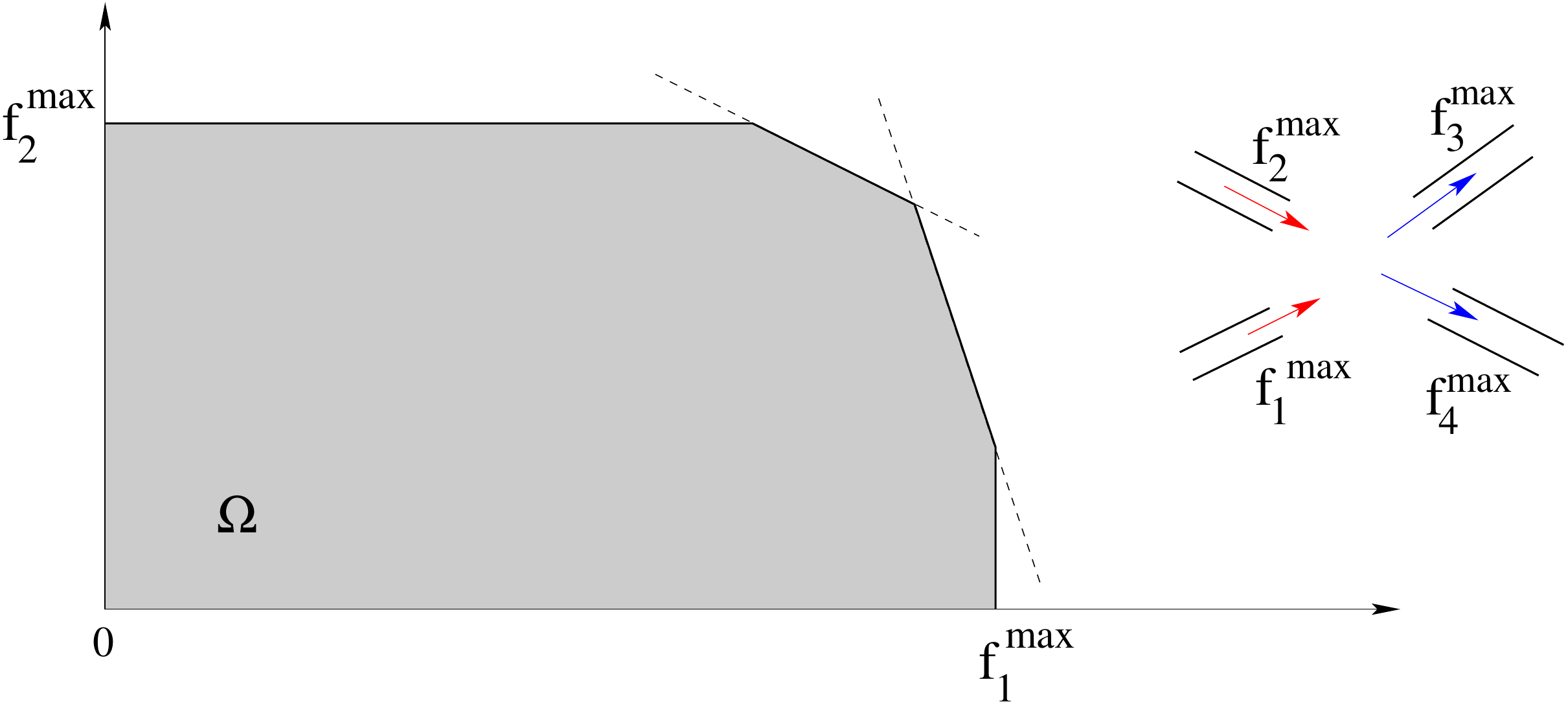}
    \caption{\small The region $\Omega\subset\R^m$ 
    of admissible incoming fluxes,  defined at (\ref{Omdef}).}
\label{f:tf63}
\end{figure}

\begin{figure}[htbp]
\centering
\includegraphics[width=0.9\textwidth]{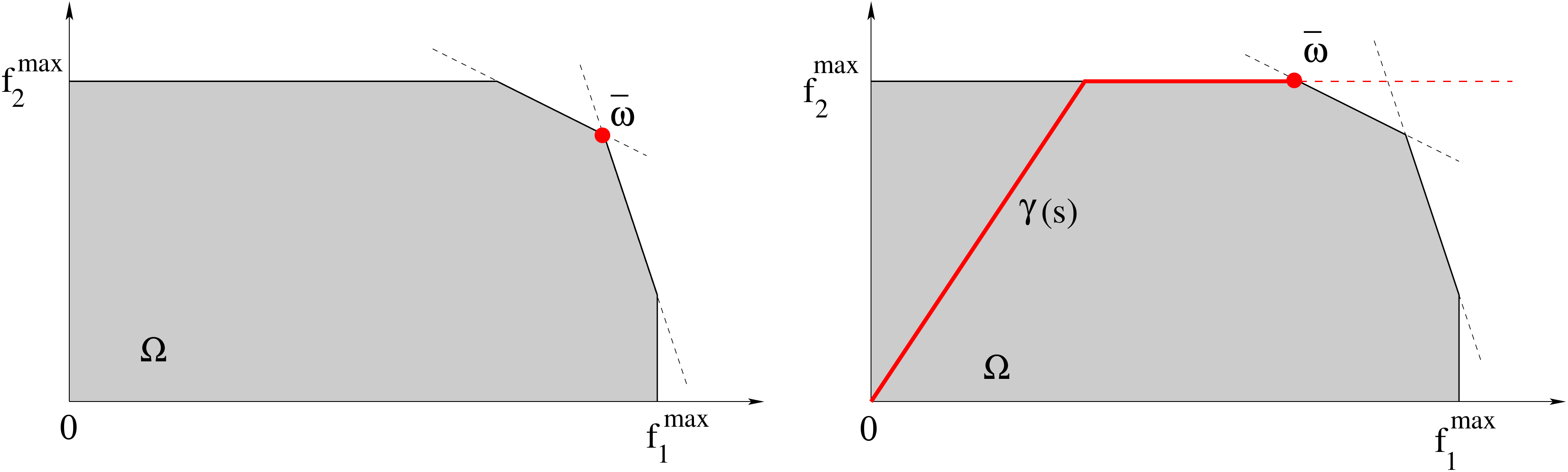}
    \caption{\small Left: the point $\bar \omega\in \Omega$ 
which    maximizes the total flux through the intersection.  
    Right: the point $\bar \omega\in \Omega$ constructed 
    by the Riemann Solver at (\ref{bbf})-(\ref{bars}).
    }
\label{f:tf240}
\end{figure}

{\bf Step 4.} To construct a 
Riemann solver, it now suffices to give a rule for selecting a point $\bar \omega=(f_1,\ldots, f_m)$
in the feasible region $\Omega$.
In general, this rule will depend on the priority coefficients 
$c_1,\ldots, c_m$ assigned to incoming roads. 
Observe that, as soon as the incoming fluxes $f_i$, $i\in\I$, are given, 
the outgoing fluxes are uniquely determined by the identities
\bel{of}
f_j~=~\sum_{i\in\I} f_i \theta_{ij}\,,\qquad\qquad j\in \O.\eeq
%For an intersection with 2 incoming and 2 outgoing roads, 
Various ways to define a Riemann Solver are illustrated by the following examples.
\v
{\bf Example 1:}  %maximizing the flux through the intersection.
Given priority coefficients $c_1,\ldots, c_m$, following \cite{CGP} 
one can choose the 
vector of incoming fluxes 
\bel{argm}\bar \omega~=~(f_1,\ldots, f_m)~\doteq~\argmax_{\omega\in \Omega}~
\sum_{i\in\I} c_i f_i\,.\eeq
In particular, if $c_1=\cdots=c_m={1\over m}$, this means we are maximizing the total flux
through the intersection (see Fig.~\ref{f:tf240}, left). 
\v
Since in (\ref{argm}) we are maximizing a linear function over a polytope,  in some 
cases
the maximum can be attained at multiple points.  This somewhat restricts the 
applicability of this model.  An alternative model, with better continuity properties, is
considered below.
\v
{\bf Example 2:} Given positive coefficients $c_1,\ldots, c_m$ as in  (\ref{tca}),
consider the one-parameter curve
$$s~\mapsto~\gamma(s) ~=~(\gamma_1(s),\ldots, \gamma_m(s)),$$
where 
$$\gamma_i(s)~\doteq~\min\{ c_i s\,,~ f_i^{max}\}.$$
As shown in Fig.~\ref{f:tf240}, right, we then choose the vector of incoming fluxes 
\bel{bbf}  \ov\omega~=~(f_1,\ldots,f_m),\qquad\qquad f_i~=~\gamma_i(\bar s),\eeq
where
\bel{bars}\bar s~=~\max~\left\{ s\geq 0\,;~~\sum_{i\in \I} \gamma_i(s)\, \theta_{ij}~\leq~
f_j^{max}
\quad\forall j\in\O\right\}.\eeq

\v

\begin{figure}[htbp]
\centering
\includegraphics[width=0.9\textwidth]{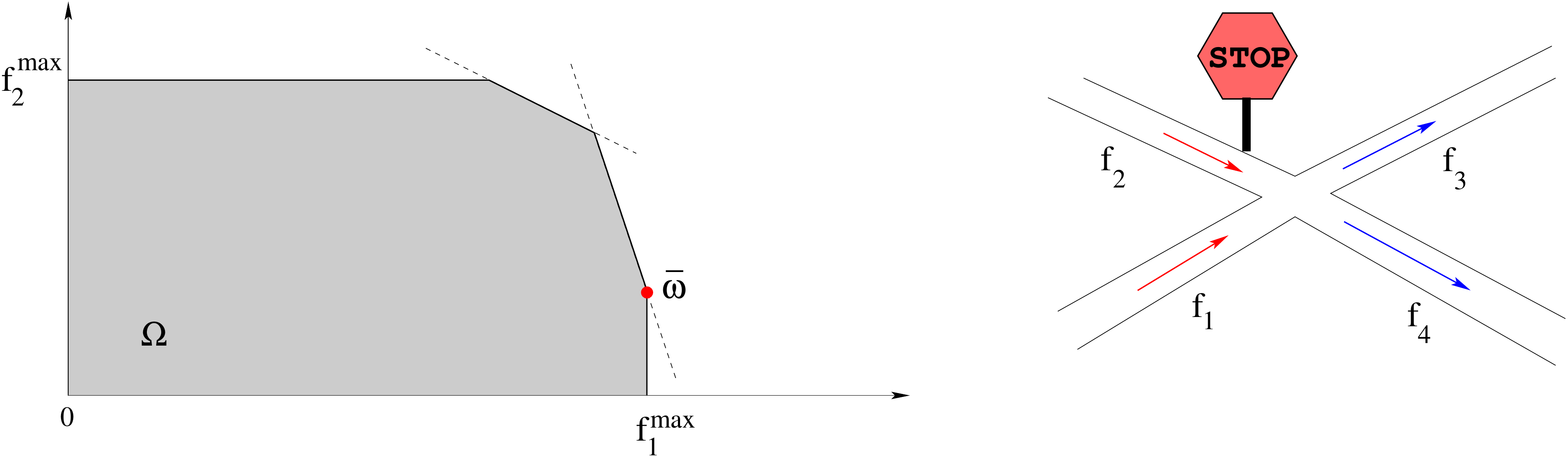}
    \caption{\small  A Riemann Solver modeling an intersection
where the second road has a stop sign.
    }
\label{f:tf212}
\end{figure}
\v
{\bf Example 3 :} To model an intersection with two incoming and two outgoing roads,
where road 2 has a stop sign, we choose the point $\ov\omega=(f_1,f_2)$
according to the following rules (see Fig.~\ref{f:tf212}).
\bel{sr1}f_1~=~\max~\{ \omega_1\,;~~(\omega_1,0)\in\Omega\}.\eeq
\bel{sr2}f_2~=~\left\{\bega{cl} 0\qquad &\hbox{if} \quad f_1<f_1^{max},\\[3mm]
\max~\{ \omega_2\,;~~(f_1,\omega_2)\in\Omega\}\qquad&\hbox{if} \quad f_1=f_1^{max}
\enda\right.\eeq
According to (\ref{sr1}),  as many cars as possible are allowed to 
arrive from road 1. According to (\ref{sr2}), if any available space is left, 
cars arriving from road 2 are allowed through the intersection.

\v
\subsection{Intersection models with buffers.}

Having defined a way to solve each Riemann problem,
a major issue is whether the Cauchy problem
with general initial data is well posed.   Assuming that the turning preferences
$\theta_{ij}$
remain constant in time, some results in this direction can be found in \cite{CGP}.

We remark, however, that in general
these turning preferences may well vary in time.   
%For example, consider a network
%as in Fig.~\ref{f:tf99}. 
%A first group of drivers originate  from road 1 eventually turn into road 4,
%while a second group of drivers
%originate from road 2 eventually turn into road 5.  In such case, the fraction 
%of cars that turn right or left at the junction $B$ at any given time $t$  will depend on how many of these
%drivers belong to the first or second group.  
%This may well change in time.
One should thus  regard  $\theta_{ij}=\theta_{ij}(t,x)$
as variables.   Assuming that drivers know in advance their itinerary,
the conservation of the number of drivers on 
road $i$ that will eventually turn into road $j$
is expressed by the additional conservation law
%the  conservation of the number of drivers  from road $i$ that will eventually turn into
%road $j$ is expressed by
\bel{cl4}[\rho_i \theta_{ij}]_t +[ \rho_i v_i (\rho)\theta_{ij}]_x~=~0.\eeq
Combining (\ref{cl4}) with the conservation law
$$(\rho_i)_t +[\rho_i v_i (\rho)]_x~=~0,$$
one obtains a linear transport equation for each of the quantities $\theta_{ij}$, namely
\bel{cl5}
(\theta_{ij})_t + v_i(\rho) \,(\theta_{ij})_x~=~0,\qquad\qquad
i\in \I, ~~j\in \O.\eeq

A surprising counterexample constructed in \cite{BY} shows that,
for a very general class of Riemann Solvers, 
one can construct measurable initial  data
$\rho_i(0,\cdot)$, $\rho_j(0,\cdot)$, and 
 $\theta_{ij}(0,\cdot)$,
 so that the Cauchy problem has two distinct entropy-admissible solutions.

The ill-posedness of these model equations represents a serious obstruction, 
toward the existence of globally 
optimal traffic patterns, or Nash equilibria, on a general network of roads.
To cope with this difficulty, in \cite{BN1} an alternative model 
was proposed, for traffic flow at an intersection. Namely, it is assumed that the junction contains a buffer (say, a traffic circle). 
Incoming cars are admitted at a rate depending of the amount of free space left in the buffer, regardless of their destination. Once they have entered the intersection, cars flow out at the maximum rate allowed by the outgoing road of their choice.    

More precisely, 
consider a constant 
$M>0$, describing the maximum number of cars that can occupy the intersection
at any given time, and constants $c_i>0$, $i\in \I$,
accounting for priorities
given to different incoming roads.  For $j\in \O$, at any time $t$ we denote by 
$q_j(t)\in [0,M]$ the number of cars, already  within the buffer, that seek to turn into road $j$.

As before, let $f_i^{max}$ and $f_j^{max}$ the maximum fluxes that can exit from 
road $i\in\I$, or can enter into road $j\in\O$.
We then require that the  incoming fluxes $f_i$ satisfy
\bel{bff}
 f_i~=~\min~\left\{ f^{max}_i\,,~~c_i\Big(M-\sum_{j\in\O} q_j\Big)\right\},
\qquad\qquad i\in\I\,.\eeq
In addition, the outgoing fluxes $ f_j$ should satisfy 
\bel{bqj}\left\{
\bega{l}\hbox{if $q_j>0$, then  $ f_j =f^{max}_j$,}\cr\cr
\hbox{if $q_j=0$, then $f_j = 
\min\Big\{ f^{max}_j, ~\sum_{i\in \I}  f_i\theta_{ij}\Big\}$,}
\enda\right.\qquad\qquad  j\in \O\,.\eeq
Having determined the incoming and outgoing fluxes $f_i$, $f_j$,
the time derivatives of  the queues $q_j$ are then computed by
\bel{dqj}
\dot q_j~=~\sum_{i\in \I} f_i \theta_{ij} - f_j\,,\qquad\qquad j\in\O.\eeq
 The well-posedness of the intersection
model with buffers, for general $\L^\infty$ data,  was proved in \cite{BN1}.

It is interesting to understand the relation between 
the intersection model with buffer, and the models based on
a Riemann Solver.   The analysis in \cite{BNo} shows that,
letting the size of the buffer $M\to 0$,
the solution of the problem with buffers converges to the solution determined by the
Riemann Solver at (\ref{bbf})-(\ref{bars}), described in Example 2.

\begin{figure}[htbp]
\centering
\includegraphics[width=0.5\textwidth]{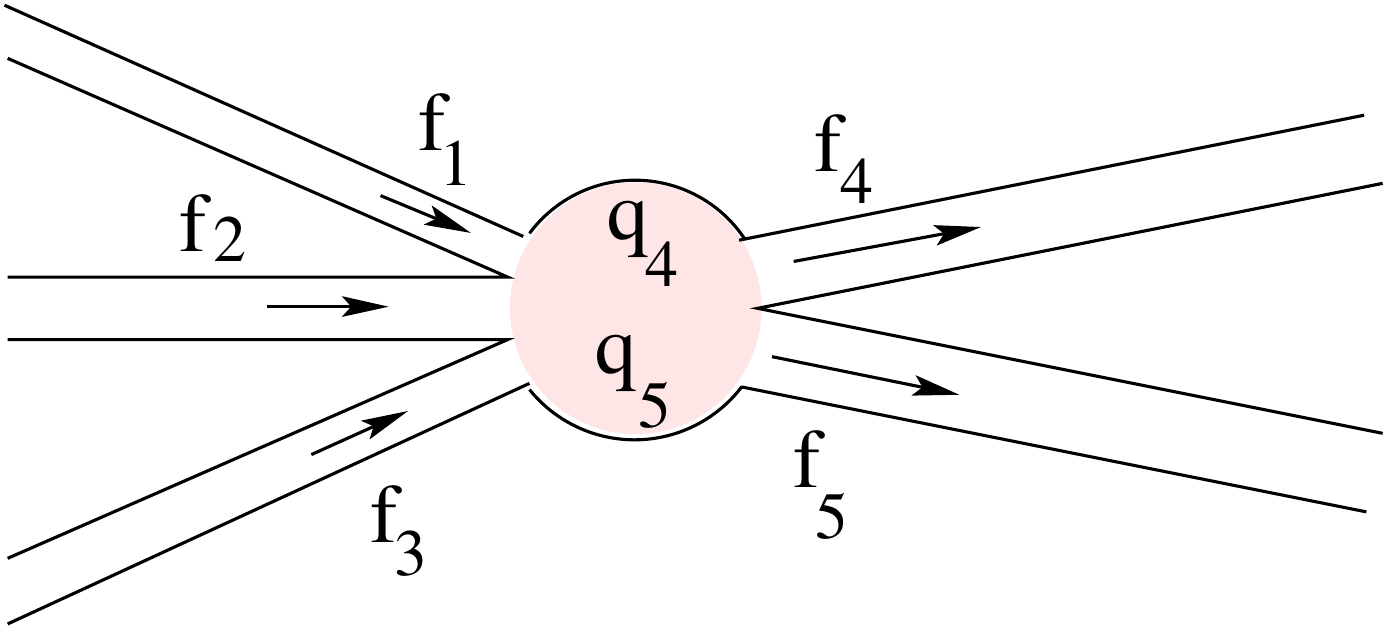}
    \caption{\small  An intersection model with a buffer. Here the queue sizes 
    $q_4, q_5$ account for the number  
 of cars that have already accessed the intersection, and are waiting to exit into roads 4 and 5 respectively. 
    }
\label{f:tf241}
\end{figure}

\subsection{Optima and equilibria on a network of roads.}
Consider a network of roads, with several intersections.   We call 
$\gamma_k$, $k=1,\ldots,\bar k$ the  arcs corresponding to the various roads, and 
$A_1,\ldots, A_\nu$ the nodes corresponding to intersections. 
It is assumed that, on the $k$-th road,  the flux function 
has the form $f_k(\rho) = \rho\, v_k(\rho)$, with $v_k$ a decreasing function of the 
density. 
As in the previous sections, traffic flow at each intersection can be modeled
in terms of a Riemann Solver, or by means of a buffer.

\begin{figure}[htbp]
\centering
\includegraphics[width=0.6\textwidth]{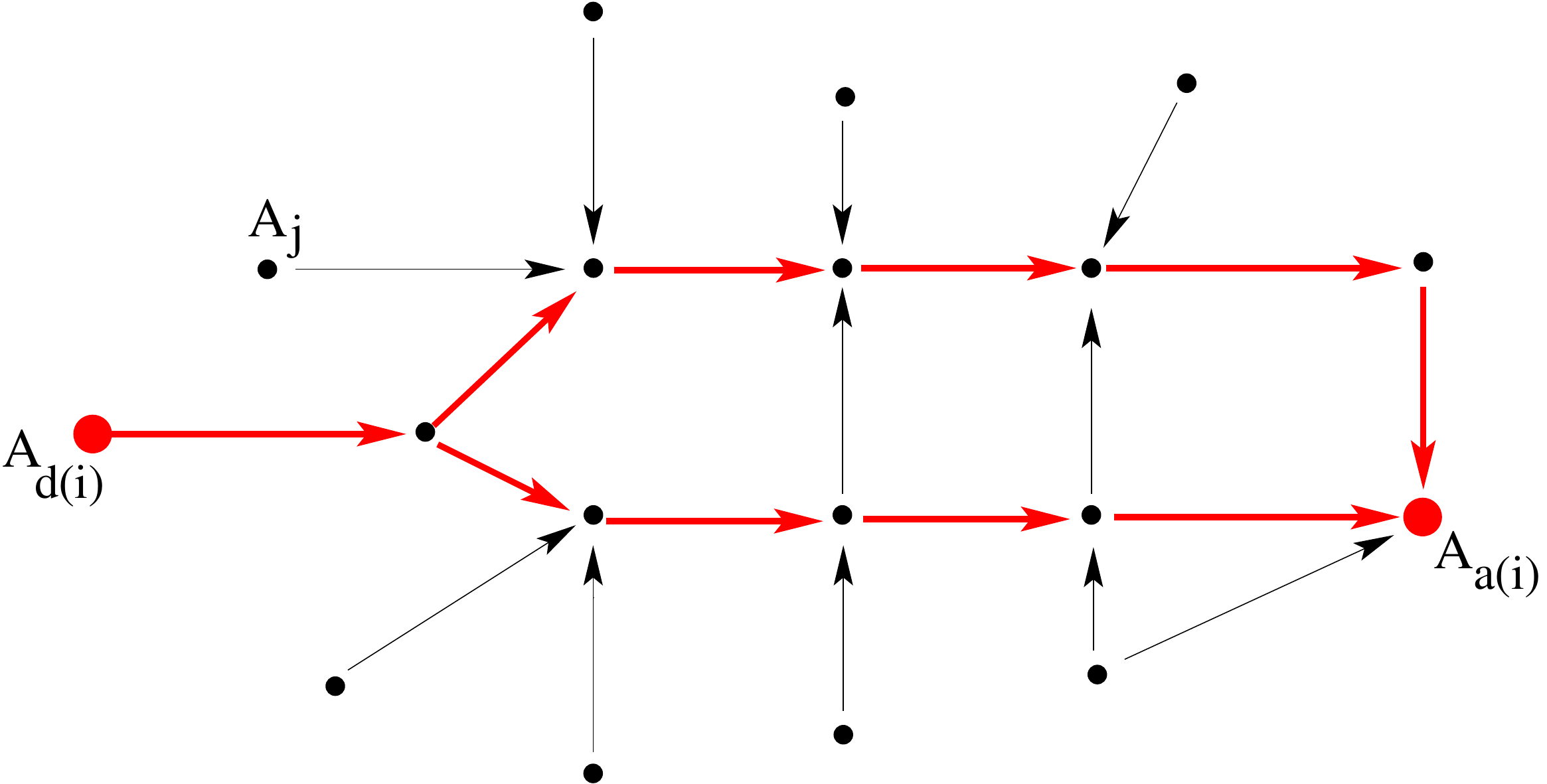}
    \caption{\small  A network of roads, with several intersections at nodes $A_j$.   
    Drivers of the $i$-th group depart from the node $A_{d(i)}$ and arrive to node $A_{a(i)}$.
  In this example  they can choose two distinct paths to reach destination. }
\label{f:tf242}
\end{figure}

We consider
$N$ groups of drivers with different origins and destinations, and possibly different 
departure and arrival costs.  As shown in Fig.~\ref{f:tf242}:
\begi
\item Drivers in the $i$-th group  depart from the node $A_{d(i)}$   and arrive at
the node $A_{a(i)}$.
\v
\item Their cost for departing at time $t$ is $\vp_i(t)$, while their arrival cost is $\psi_i(t)$.
\item They can use different paths $\Gamma_1, \Gamma_2,\ldots$  to reach destination.
\endi

In the following, 
$\bar u_{i,p}(\cdot)$ will denote the departure rate of drivers of the $i$-th group, 
who choose the path $\Gamma_p$  to reach destination.
Calling $G_i$ the total number of drivers in the $i$-th group,
we say that the departure rates $\bar u_{i,p}$ are {\bf admissible} if, for every $i=1,\ldots,N$
they satisfy the obvious constraints
\bel{bua}\bar u_{i,p}(t)~\geq~0,\qquad\qquad
\sum_p\int_{-\infty}^{+\infty} \bar u_{i,p}(t)\, dt~=~G_i\,.\eeq
Given the departure rates, in principle one can then solve the equation of traffic flow on the 
whole network and determine the arrival times of the various drivers.
We call
$$\tau_p(t)~=~\hbox{arrival time of a driver departing at time $t$, traveling along the 
path $\Gamma_p$.}$$
With these notations, we can introduce

\begin{definition} 
 An admissible family
$\{ \bar u_{i,p}\}$ of departure rates is {\bf globally optimal}
if it minimizes the sum of the total costs of all drivers
$$J(\bar u) ~\doteq~\sum_{i,p}\int \Big( \vp_i(t) + \psi_i(\tau_{p}(t))\Big) 
\bar u_{i,p}(t)\, dt\,.$$
\end{definition}

\begin{definition} 
An admissible family
$\{ \bar u_{k,p}\}$ of departure rates is a 
{\bf Nash equilibrium}
if no driver of any group can lower his own total  cost by 
changing departure time or switching to a different path to reach destination.
\end{definition}

From the above definition it follows the existence of constants $C_1,\ldots, C_N$
such that
$$\bega{l}\vp_k(t) + \psi_k(\tau_{p}(t))~=~C_k\qquad\qquad \forall t\in 
{\rm Supp} (\bar u_{k,p})\,,\cr
\cr
\vp_k(t) + \psi_k(\tau_{p}(t))~\geq ~C_k\qquad\qquad \forall t\in \R\,.\enda$$
As remarked in the Introduction, a similar characterization for 
the globally optimal solution is much harder to justify.

\v

\begin{figure}[htbp]
\centering
  \includegraphics[scale=0.45]{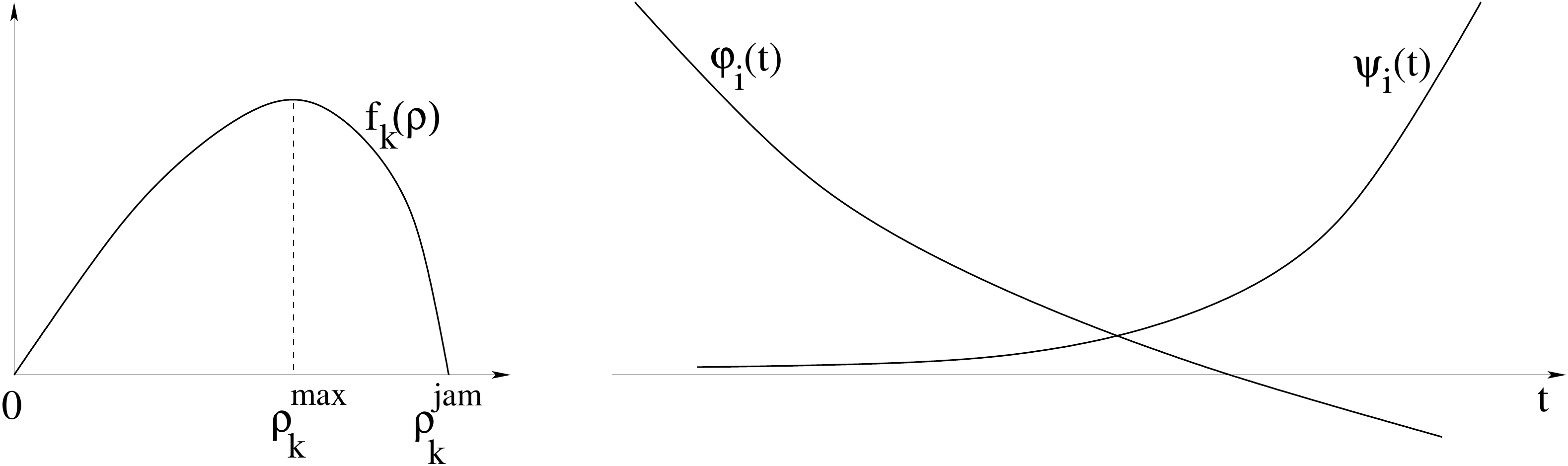}
    \caption{\small A flux function $f_k$,  a departure  cost 
    function $\vp_i$, and and arrival cost function $\psi_i$, 
    satisfying the assumptions (A1)-(A2).}
\label{f:tf243}
\end{figure}

In the above setting, a natural set of assumptions is (see Fig.~\ref{f:tf243})
\begin{itemize}
\item[{\bf (A1)}]  On each road $k=1,\ldots,\bar k$, the flux function $f_k$ satisfies
\bel{fi2}
f_k\in \C^2, \qquad f_k''<0,\qquad f_k(0)= f_k(\rho_k^{jam}) = 0.\eeq

\item[{\bf (A2)}] For each group of drivers $i=1,\ldots,N$, the cost functions $\vp_i$, $\psi_i$ satisfy
\bel{vpi}\vp_i'~<~0,\qquad \psi_i, \psi_i'~> ~0,\qquad\quad
 \lim_{|t|\to \infty} \Big(\vp_i(t)+ \psi_i(t)\Big)~=~+\infty.\eeq
\end{itemize}

When all intersections are modeled in terms of a buffer as in (\ref{bff})--(\ref{dqj}),
under the assumptions {\bf (A1), (A2)},
 the existence of
at least 
one globally optimal solution and at least one Nash equilibrium solution
was proved in \cite{BN2}.

\v

\section{Optimal solutions: a single road, several groups of drivers}
\label{sec:3}
\setcounter{equation}{0}

Consider a single road, where the traffic density is governed by the conservation law
\bel{cl1}
\rho_t + f(\rho)_x~=~0\qquad\qquad \text{for}\quad x\in [0,L].\eeq
We assume that $N$ groups of drivers are present, of sizes $G_1,\ldots, G_N$, with 
departure and arrival costs $\vp_i$, $\psi_i$, $i=1,\ldots,N$.
The flux function will be denoted by
$$u(t,x)~=~f(\rho(t,x))~=~\sum_{i=1}^N u_i(t,x).$$
Here 
\bel{uidef}u_i(t,x)~=~ \theta_i(t,x)\, u(t,x)\eeq
is the flux of drivers of the $i$-th group.
As in (\ref{tca}), we always assume that
\bel{tid}
\theta_i(t,x)\,\geq\,0,\qquad\qquad\sum_{i} \theta_i(t,x)\,=\,1.\eeq
For each $i\in \{1,\ldots,N\}$, the conservation of the 
number of drivers of the $i$-th family yields 
the additional conservation law
\bel{2}
(\rho \theta_{i})_t +\Big( \rho v(\rho)\,\theta_{i}\Big)_x~=~0\,.\eeq
By (\ref{cl1}), one obtains the  linear equations
\bel{3}
\theta_{i,t} + v(\rho) \,\theta_{i,x}~=~0,\qquad\qquad i=1,\ldots,N.\eeq

The incoming flux at the beginning of the road is
\bel{bc1}u(t,0)~=~\bar u (t)~=~\sum_{i=1}^N \bar \theta_i(t)\, \bar u(t).\eeq
The global optimization problem can be formulated as follows.

\begi
\item[{\bf (OP)}] Given the constants $G_i>0$, $i=1,\ldots,N$, find departure rates
$\bar u_i(t)=\bar\theta_i(t)\bar u(t)$ which provide an optimal solution to the problem
\bel{TC}\hbox{minimize:}\qquad J(\bar u_1,\ldots,\bar u_N)~\doteq~
\sum_{i=1}^N~\int_{-\infty}^{+\infty}\Big[ u_i(t,0) \vp_i(t) + u_i(t,L)\psi_i(t)\Big]\, dt ,\eeq
and satisfy the constraints
\bel{TConst} \bar u(t)\in [0,M],\qquad
\bar \theta_i(t)\,\geq \,0,\qquad\sum_{i=1}^N \bar \theta_i(t)~=~1,
\qquad\forall t\in\R,\eeq
\bel{CK}\bar u_i(t)\,\geq \, 0,\qquad\qquad \int_{-\infty}^{+\infty} \bar u_i(t)\, dt~=~
G_i\,,\qquad \quad i=1,\ldots,N.\eeq
\endi

We recall that $u_i(t,L)=\theta_i(t,L) u(t,L)$ is the rate at which the drivers of the $i$-th group
arrive at the end of the road.

Since no intersections are present, the existence of a globally optimal solution
follows as  a special case of the result in \cite{BN2}.  Here we briefly recall the 
main argument in the proof.

{\bf 1.} Let $(\bar u_1^{(n)},\ldots, \bar u_N^{(n)})_{n\geq 1}$ be a  minimizing sequence
of admissible departure rates. Namely
$$\bar u^{(n)}_{i}(t)\geq 0, \qquad\qquad
 \sum_{i=1}^N\bar u^{(n)}_{i}(t)~\leq M,\qquad\qquad
\int_{-\infty}^{+\infty} \bar u^{(n)}_{i}(t)\, dt~=~G_i$$
for every  $n\geq 1$,  and moreover
$$\lim_{n\to\infty}~ J(\bar u_1^{(n)},\ldots,\bar u_N^{(n)})~=~
\inf J(\bar u_1,\ldots,\bar u_N).$$
{\bf 2.} By the assumption {\bf (A2)}, as $t\to \pm\infty$ the cost functions 
$\vp_i,\psi_i$ become very large.
By possibly modifying the functions $\bar u_i$, we can thus 
obtain a minimizing sequence
where all departure rates vanish outside a fixed time interval $[a,b]$.

{\bf 3.} By taking a subsequence, we obtain a weak limit 
$ (\bar u_1^{(n)},\ldots, \bar u_N^{(n)})\wto (\bar u_1,\ldots, \bar u_N)$.

The boundedness of the supports guarantees that these limit departure rates are
still admissible (i.e., no mass leaks at infinity).  

{\bf 4.} Call $u_i^{(n)}(t,x)$, $u_i(t,x)$ the corresponding solutions.
By the genuine nonlinearity of the conservation law (\ref{cl1}), after taking a subsequence,
  one obtains the strong
convergence $u^{(n)}(\cdot, L)\to u(\cdot,L)$ in $\L^1 (\R)$,
and the weak convergence of the departure and arrival rates 
$$u^{(n)}_i(\cdot,0)\wto u_i(\cdot,0),\qquad\qquad 
u^{(n)}_i(\cdot,L)\wto u_i(\cdot,L).$$  Since the cost functional in (\ref{TC})
is linear w.r.t.~these departure and arrival rates, it is continuous
w.r.t.~weak convergence. This yields the optimality
of the departure rates $(\bar u_1,\ldots, \bar u_N)$.

\subsection{Optimality conditions.}
In the remainder of this section, 
we seek necessary conditions for a solution to be optimal.
As a first step, we derive an explicit representation of the solution.

\begin{figure}[htbp]
   \centering
 \includegraphics[width=1.0\textwidth]{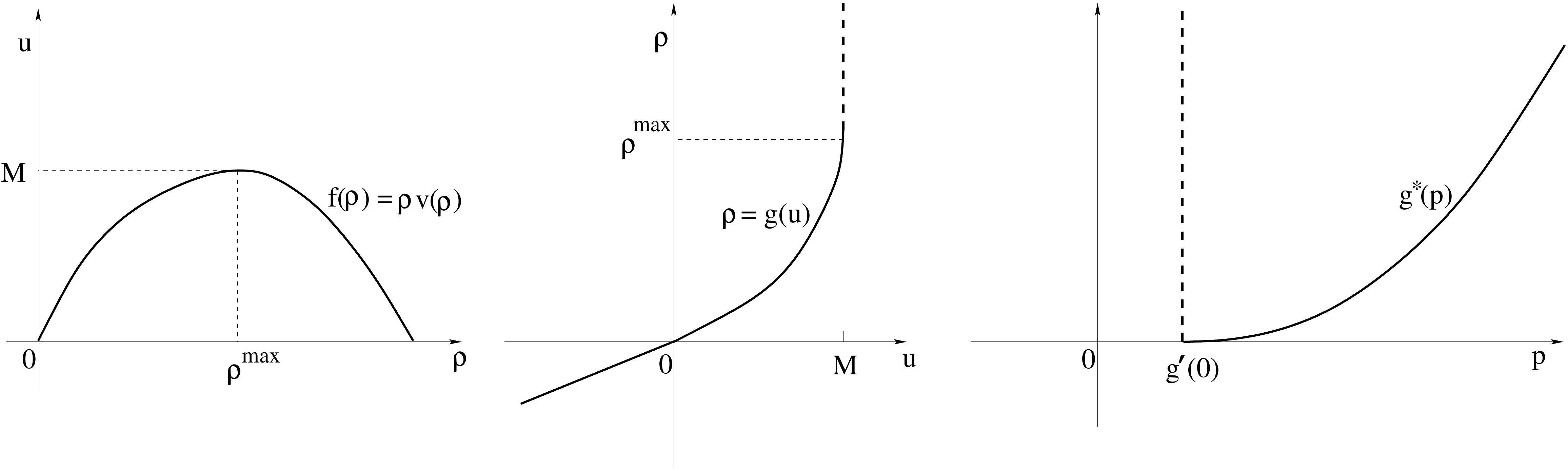}
   \caption{\small  Left: the function $\rho\mapsto
   \rho\,v(\rho)$ describing the flux of cars. Middle:
the  function $g$, implicitly defined
by $g(\rho v(\rho)) = \rho$ and extended according to (\ref{gext}).
Right: the Legendre transform
   $g^*$.}
   \label{f:tf244}
\end{figure}

Following \cite{BH1, BH2}, it is convenient to switch the roles of the variables $t,x$, and
write the density $\rho$ as a function of the flux $u$.
The boundary value problem (\ref{cl1})-(\ref{bc1}) thus becomes a
Cauchy problem
for the conservation law describing the flux $u=\rho v(\rho)$,
namely
\bel{CP}
u_x + g(u)_t ~=~0\,,\eeq
\bel{inc}
u(t,0)~=~\bar u(t)\,.\eeq
As shown in Fig.~\ref{f:tf244},
the function
$u\mapsto g(u)~=~\rho$
is defined as a partial 
inverse of the function $\rho\mapsto \rho\, v(\rho) =u$,
assuming that  
$$0~\leq~u~\leq~ M~\doteq~\max_{\rho\geq 0}  \rho v(\rho)\,,
\qquad\qquad  0~\leq~\rho~\leq ~\rho^{max}.$$
For convenience, we extend $g$ to the entire real line by setting
\bel{gext}
g(u)~\doteq~\left\{ \bega{cl} g'(0+) u\qquad &\hbox{if}~~u<0,\\[3mm]
+\infty\qquad &\hbox{if}~~u> M\,.\enda\right.\eeq

The solution to (\ref{CP})-(\ref{inc}) can now be expressed by means of the Lax formula
\cite{Evans, Lax}.
Namely, call 
\bel{Leg}
g^*(p)~\doteq~\max_{u}  \{ pu - g(u)\}\eeq
the Legendre transform of $g$.   Notice that 
$$g^*(p)~=~+\infty\qquad\hbox{for}\quad p< g'(0).$$
On the other hand, for $p \geq g'(0)$ the strict convexity of $g$ implies that there exists
a unique value $u= \gamma(p)\geq 0$ where the maximum in (\ref{Leg}) is attained, 
so that
$$g^*(p)~=~p\cdot \gamma(p) - g(\gamma(p)).$$
This function $\gamma: [g'(0), +\infty[\,\mapsto [0, M[\,$ is implicitly defined by the relation
\bel{gg}
g'(\gamma(p)) ~=~p.\eeq
Consider the integrated function
$$U(t,x)~\doteq~\int_{-\infty}^t u(\tau,x)\, d\tau,$$
which measures the number of drivers that  have crossed the point $x$
along the road before time $t$.  
The conservation law (\ref{CP}) can be equivalently written as 
a Hamilton-Jacobi equation
\bel{HJ}U_x + g(U_t)~=~0\eeq
with data at $x=0$
\bel{HJ0}U(t,0)~=~\ov U(t)~=~\int_{-\infty}^t \bar u(s)\, ds.\eeq
The solution to (\ref{CP})-(\ref{inc}) is now provided by the Lax formula
\bel{LF1}
U(t,x)~=~\min_\tau ~\left\{ x g^*\Big( {t-\tau\over x}\Big) + \ov U(\tau)\right\}\,,\eeq
\bel{LF2} \tau(t,x)~\doteq~\argmin_\tau 
~\left\{ x g^*\Big( {t-\tau\over x}\Big) + \ov U(\tau)\right\}\,,\eeq
\bel{LF3}
u(t,x)~=~\gamma\left( t-\tau(t,x)\over x\right).\eeq
We observe that the function $U=U(t,x)$ is globally Lipschitz continuous. Its values
satisfy
\bel{LF4}
U(t,x)~\in~[0,G],\qquad\qquad G=G_1+\ldots+G_N\,.\eeq

Car trajectories $t\mapsto y(t)$ are defined to be the solutions to
the ODE
\bel{ctra} \dot y(t)~=~v(\rho(t, y(t))).\eeq
In the region where $\rho>0$, and hence $u=\rho v(\rho)>0$ as well, these trajectories
coincide with the level curves of the integral function $U$.
Indeed, observing that $v = u/\rho$, when $\rho = g(u)>0$ we can write
$$\dot y(t)~=~v(\rho(t, y(t)))~=~{u(t, y(t))\over g(u(t, y(t))}~=~{U_t(t, y(t))\over g(U_t(t, y(t))}
\,.$$
By (\ref{HJ}) one has
\bel{CT}{d\over dt} \,U(t, y(t))~=~U_t + U_x \dot y(t)~=~U_t + U_x \,{U_t\over g(U_t)}~=~0.
\eeq
More generally, consider a car departing at time $t_0$.   The solution to the Cauchy problem
\bel{ct3} \dot y(t)~=~v(\rho(t, y(t))),\qquad\quad y(t_0)\,=\,0\eeq
can be determined by the formula
\bel{ct4}
y(t)~=~\inf\left\{ x\,;~~U(t,x)>\ov U(t_0)\quad\hbox{or}\quad x={t-t_0\over v(0)}\right\}.\eeq
The arrival time of a driver departing at time $t_0$ is
\bel{taf}\tau^a(t_0)~=~\sup\left\{ t'\,;~~U(t',L)<\ov U(t_0)\quad\hbox{or} ~~t'= t_0 + {L\over v(0)}
\right\}.\eeq
By (\ref{3}), the functions $\theta_i$ are constant along car trajectories.

In order to compute the arrival rates $u_i=\theta_iu$ at the terminal 
point of the road $x=L$, we 
first observe that the map 
$$t~\mapsto~\ov U(t)$$ in  (\ref{HJ0}) is nondecreasing. Hence we can define an inverse
by setting
\bel{Thi}
\ov \tau(s)~\doteq~\inf~\{t\in\R\,;~\ov U(t)\geq s\}~\in~\Big[0\,,~ \sum_i G_i\Big].\eeq
We then introduce the functions
\bel{THI}
\Theta_i(s)~\doteq~\bar\theta_i(\bar\tau(s))
\eeq
By (\ref{3}), the functions $\theta_i=\theta_i(t,x)$ are constant along car trajectories.
The general solution to  (\ref{3})-(\ref{bc1}) can thus be written as
\bel{TIS}
\theta_i(t,x)~=~\Theta_i(U(t,x)).\eeq

%\int_{-\infty}^{\ov\tau(s)} \bar\theta_i(t)\, \ov U(t)\, dt.\eeq
%\bel{TH'}
%\Theta_i'(s)~=~{\partial\over\partial s} \Theta_i(s)~=~\bar \theta_i(\ov\tau(s)).\eeq
%Notice that, among the first $s$ drivers who depart, $\Theta_i(s)$ gives the amount of those who belong to the $i$-th group.

We shall be mostly interested in the terminal values $u(t,L)$, describing the rate 
at which cars arrive at the end of the road.
%To shorten the notation, the initial point on the characteristic through 
%$(t,L)$ will be denoted by
%$(\tau(t),0)$, where
%\bel{TT}\tau(t)~\doteq~\argmin_\tau ~\left\{ x g^*\Big( {t-\tau\over L}\Big) + \ov U(\tau)\right\}\eeq
Denoting the arrival distribution as $U^a (t)\doteq U(t,L)$, 
the total cost can now be written as
\bel{cost}
J~=~\int\vp(t) \,d \ov U(t)+ \sum_i \int \psi_i(t)\, \Theta_i(U^a(t)) \,d U^a(t).\eeq

\begin{theorem}\label{t:1}
 Let the flux function $f$ satisfy the standard assumptions (\ref{fi}).
 Assume that all drivers have the same departure cost $\vp_i=\vp$
and possibly different arrival costs $\psi_1,\ldots,\psi_N$, satisfying (\ref{vpi}).
Let $(\bar u_1,\ldots, \bar u_N)$ an optimal departure rate,
minimizing the total cost to all drivers.  

Then  the corresponding solution does not contain shocks.
Moreover, there exists constants $C_1,\ldots, C_N$
such that, setting  
\bel{pps}
\psi(t)~\doteq~\min_k ~(\psi_k(t)-C_k)\,,\eeq
the following holds:  
\begi
\item[(i)]
For any $t\in \R$, let
$T(t)$ be the unique time such that
\bel{ncc}
\vp(t)+\psi(T(t))~=~0\,.\eeq
Then, for every point $(t', x')$
along the segment with endpoints $(t,0)$ and $(T(t),L)$, one has
\bel{uc}u(t',x')~=~\left\{\bega{cl} \ds\gamma\left( {T(t)-t\over L}\right)\qquad &\hbox{if}
\quad \ds {T(t)-t\over L}\,\geq \,g'(0),\\[4mm]
0\qquad &\hbox{otherwise}.\enda \right.\eeq
\item[(ii)] Calling $u_i(\cdot, L)$ the arrival rate of drivers of the $i$-th group,
one has
\bel{SUI}
\hbox{Supp}(u_i(\cdot,L))~\subseteq~\{ s\,;~~\psi(s)\,=\,\psi_i(s)-C_i\}.\eeq
\endi
\end{theorem}

%The constant $C_i$ can here be interpreted as the marginal cost for
%inserting one more driver of the $i$-th group.
A proof of Theorem~\ref{t:1} will be given in the next section.

\section{Proof of the necessary conditions}
\label{sec:4}
\setcounter{equation}{0}
As a preliminary, we review the basic theory of scalar 
conservation laws with convex flux \cite{Bbook, Evans, Smoller}.
Notice that in (\ref{CP}) the usual role of the variables $t,x$ is reversed, 
because of the particular meaning of the equations.

\begin{figure}[htbp]
   \centering
 \includegraphics[width=0.8\textwidth]{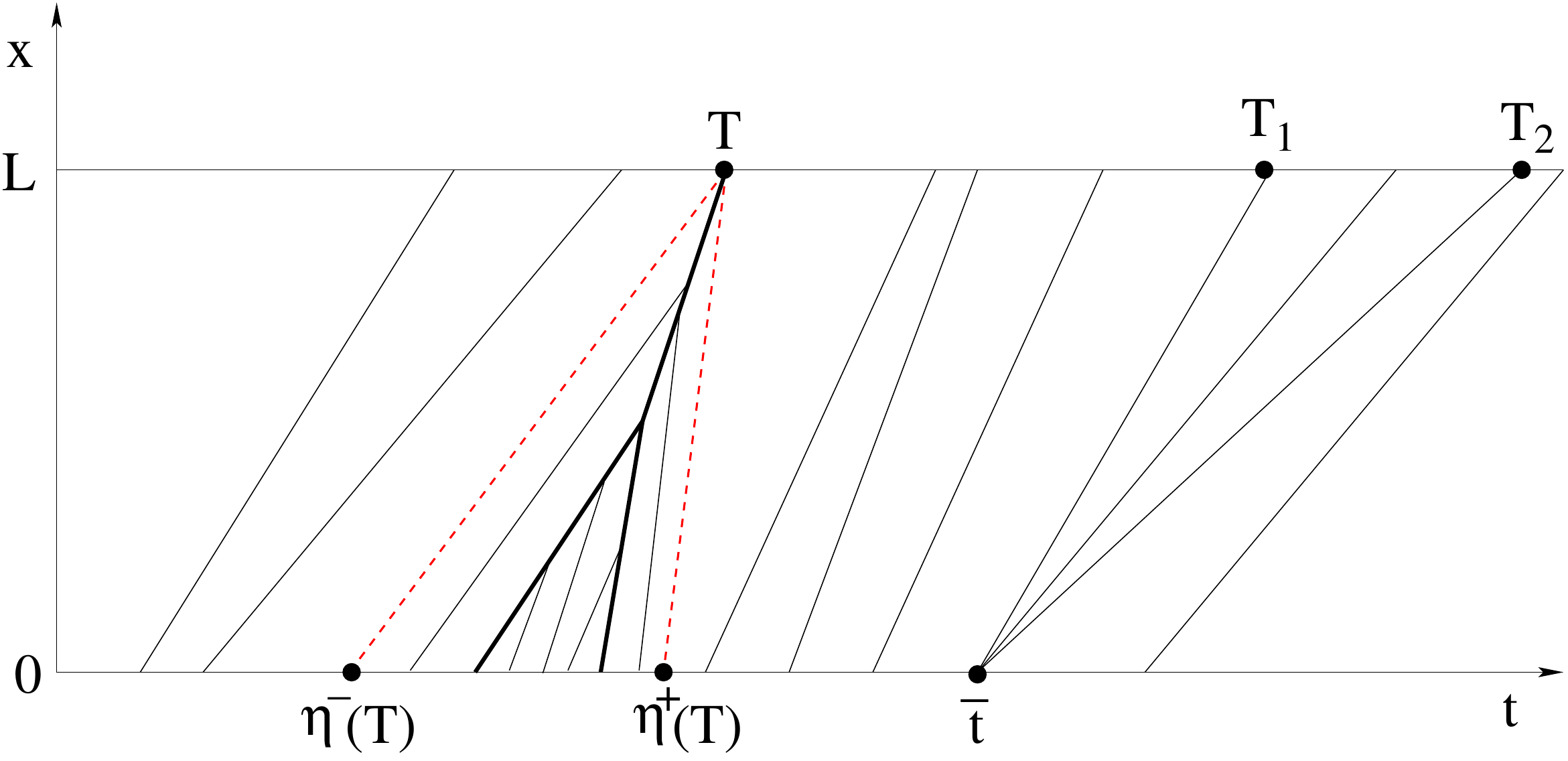}
   \caption{\small  Characteristic lines for a solution of (\ref{CP}). Here the point $(T,L)$
   lies along a shock, and all (generalized) characteristics starting 
   at a point $(t,0)$ with $t\in [\eta^-(T), \,\eta^+(T)]$ eventually reach $(T,L)$.
  Notice that there are several characteristics starting from 
  the point $(\bar t, 0)$, which is the center of a rarefaction wave.
   }
   \label{f:tf248}
\end{figure}

Let $u=u(t,x)$ be a weak solution to (\ref{CP}), taking values within the interval
$[0,M]$. This solution is entropy admissible
if it contains only downward jumps, namely
$$u(t+,x)~\doteq ~\lim_{s\to t+} u(s,x)~\leq~\lim_{s\to t-} u(s,x)~\doteq~u(t-, x).$$
By a generalized characteristic we mean a function 
$x\mapsto t(x)$ which provides a solution
to the differential inclusion
\bel{DI} {d\over dx} t(x)~\in ~\Big[ g'(u(t+,x)),\, g'(u(t-,x))\Big].\eeq
For any given point $(T,L)$, there exists a minimal and a maximal 
backward characteristic.   As shown in 
Fig.~\ref{f:tf248}, we denote by $(\eta^-(T),0)$ and $(\eta^+(T),0)$ 
the initial points of these characteristics.  Calling $\ov U$ the integral function
in (\ref{LF1}), the points $\eta^-(T)$ and $\eta^+(T)$ are respectively 
the minimum and the maximum
elements within the set 
\bel{argmin}I(T)~\doteq~ \left\{t\in \R; t = \argmin\limits_{\tau} \left\{ L g^*\left({T-\tau\over L}\right) + \ov U(\tau) \right\}  \right\}\eeq
where the function $\Lambda(t)~\doteq~L g^*\left({T-t\over L}\right) + \ov U(t)$ 
attains its global minimum.

Two cases can occur:
\begi
\item[(i)] The global minimum in (\ref{argmin}) 
is attained at a single point $t^* =\eta^-(T)=\eta^+(T)$.   

The function $u(\cdot, L)$ is 
then continuous at the point $T$, and 
$$u(T,L)~=~\gamma\left( {T- t^*\over L}\right).$$

\item[(ii)] The global minimum in (\ref{argmin}) 
is attained at multiple points, hence $\eta^-(T)<\eta^+(T)$.   

 In this case the solution contains a shock
through the point $(T,L)$.   Recalling (\ref{gg}), the left and right values across the shock are
determined by
$$u(T-,L)~=~\gamma\left( {T- \eta^-(T)\over L}\right),
\qquad\qquad u(T+,L)~=~\gamma\left( {T- \eta^+(T)\over L}\right).$$
\endi
We observe that characteristics do not cross each other.  Indeed, one has the implication
\bel{ccc}
T_1\,<\,T_2\qquad\implies\qquad \eta^+(T_1)\,\leq \,\eta^-(T_2).\eeq
{}From (\ref{ccc}) it immediately follows that
\begi
\item[(i)] The profile  $u(\cdot, L)$ can contain at most countably many shocks.
Namely, there can be at most countably many points $T_i$ such that
$ \eta^-(T_i)<\eta^+(T_i)$.

\item[(ii)] There can be at most countably many points $t_i$  
such that
 \bel{fc}
 t_i~= ~\eta^+(T_1)~=~\eta^-(T_2)\,,\eeq
 for two distinct points $T_1<T_2$. 
 \endi

We recall that, given any function $\phi\in \L^1_{loc}(\R)$, almost every point $t\in\R$
is a Lebesgue point of $\phi$. By definition, this means
\bel{Leb}
\lim_{h\to 0+} {1\over h} \int_{t-h}^{t+h} |\phi(s)-\phi(t)|\, ds~=~0.\eeq

As proved in \cite{BH1}, if $t$ is a Lebesgue point of the initial datum $\bar u(\cdot)$, 
then there exists a unique forward characteristic starting at $t$.
In particular, $t$ cannot be the center of a rarefaction wave, and
there exists a unique point $T$ such that
\bel{tin} t~\in~ [\eta^-(T), \, \eta^+(T)].\eeq

The next lemma is concerned with the stability of the map $T\mapsto \eta^\pm(T)$,
w.r.t.~small perturbations in the initial datum $\bar u$.

\begin{lemma}\label{l:5}
Let $u=u(t,x)$ be the unique  entropy weak solution of (\ref{CP})-(\ref{inc}).
Assume that $t$ is a Lebesgue point for the initial datum $\bar u$, and 
let $T$ be the unique point such that (\ref{tin}) holds.
Then, for any $\ve>0$, one can find $\delta,\delta'>0$ such that the following holds.

Let $\bar u^\dagger$ be a second initial datum, with 
\bel{uud} \|\bar u^\dagger - \bar u\|_{\L^1}~\leq~\delta\,.\eeq
If we call $u^\dagger$ the corresponding solution, and define the maps 
$(\eta^\dagger)^\pm$
accordingly,
then %for any $T'\in [T-\delta', T+\delta'] $ one has
\bel{close}[t-\delta',\, t+\delta']~\subseteq~
[(\eta^\dagger)^+(T-\ve) ,\, (\eta^\dagger)^-(T+\ve) ].\eeq
 \end{lemma}

 {\bf Proof.} {\bf 1.}  Let $\ve>0$ be given. By the uniqueness assumption, for the solution
 $u$ the backward characteristics through the points $T-\ve$ and $T+\ve$
 satisfy
 $$
 \eta^+(T-\ve)~<~t~<~\eta^-(T+\ve).$$
 Hence we can find  $\delta'>0$ such that 
 \bel{bk1}
 \eta^+(T-\ve)~<~t-2\delta'~<~t+2\delta'~<~\eta^-(T+\ve).\eeq
\v
{\bf 2.} If the conclusion of the lemma does not hold, we could find a sequence 
of initial data $\bar u_n$
with $\|\bar u_n-\bar u\|_{\L^1}\to 0$, such that the corresponding maps $\eta_n^\pm$ satisfy
\bel{bk2}
\eta_n^+(T-\ve)~\geq ~t-\delta'\qquad\hbox{or}\qquad \eta_n^-(T-\ve)~\leq ~t+\delta'.\eeq
To fix the ideas, assume that the first case holds. Namely,  for every $n\geq 1$, there exists $t_n\geq t-\delta'$ such that
\bel{bk3}  L g^*\left({T-\ve -t_n\over L}\right) + \ov U_n(t_n)~=~
\min_\tau  \left\{L g^*\left({T-\ve -\tau\over L}\right) + \ov U_n(\tau)\right\}.\eeq
By possibly taking a subsequence we can assume $t_n\to t^*\geq t-\delta'$.
The uniform convergence $\ov U_n\to \ov U$ now yields
\bel{bk31}  L g^*\left({T-\ve -t^*\over L}\right) + \ov U(t^*)~=~
\min_\tau  \left\{L g^*\left({T-\ve -\tau\over L}\right) + \ov U(\tau)\right\}.\eeq
This implies $\eta^+(T-\ve)\geq t-\delta'$, reaching a contradiction.
 \endproof
  
\begin{remark}\label{r:shock}
{\rm As shown in Fig.~\ref{f:tf248}, consider a solution $u=u(t,x)$ containing a shock
through the point $(T,L)$.    Then we can modify the initial data at $x=0$ 
inside the interval $[\eta^-(T), \eta^+(T)]$ so that the solution perturbed solution $u^\dagger$
contains a centered compression wave which breaks exactly at $(T,L)$.  
In view of (\ref{gg}),
This is achieved by taking
$$\bar u^\dagger(t)~=~
\left\{\bega{cl}  \gamma\left({T-t\over L}\right) \qquad &\hbox{if}\quad
 t\in [\eta^-(T), \eta^+(T)],\\[4mm]
\bar u(t)\qquad&\hbox{if}\quad t\notin [\eta^-(T), \eta^+(T)].\enda\right.
$$
This ensures that all characteristics starting at a point $(t,0)$ with $t\in [\eta^-(T), \eta^+(T)]$
join together at the point $(T,L)$.  

Consider the constant
$$\lambda~\doteq~ L g^*\left({T-\eta^-(T)\over L}\right) + \ov U(\eta^-(T))~=~
L g^*\left({T-\eta^+(T)\over L}\right) + \ov U(\eta^+(T)).$$
Then  the corresponding integrated function $\ov U^\dagger$
satisfies  
$$L g^*\left({T-t\over L}\right) + \ov U^\dagger(t)~=~\lambda\qquad\qquad\forall 
t\in [\eta^-(T), \eta^+(T)].$$
By \eqref{argmin}, this  implies
\bel{UUD}
\ov U^\dagger(t)~\leq~\ov U(t)\qquad\qquad\forall t\in\R,\eeq
while the corresponding solutions coincide at $x=L$, namely
\bel{UU2}
U^\dagger(t,L)~=~U(t,L),\qquad\qquad u^\dagger(t,L)~=~u(t,L),\qquad\forall t\in\R.\eeq
 }
\end{remark}

\begin{figure}[htbp]
   \centering
 \includegraphics[width=0.8\textwidth]{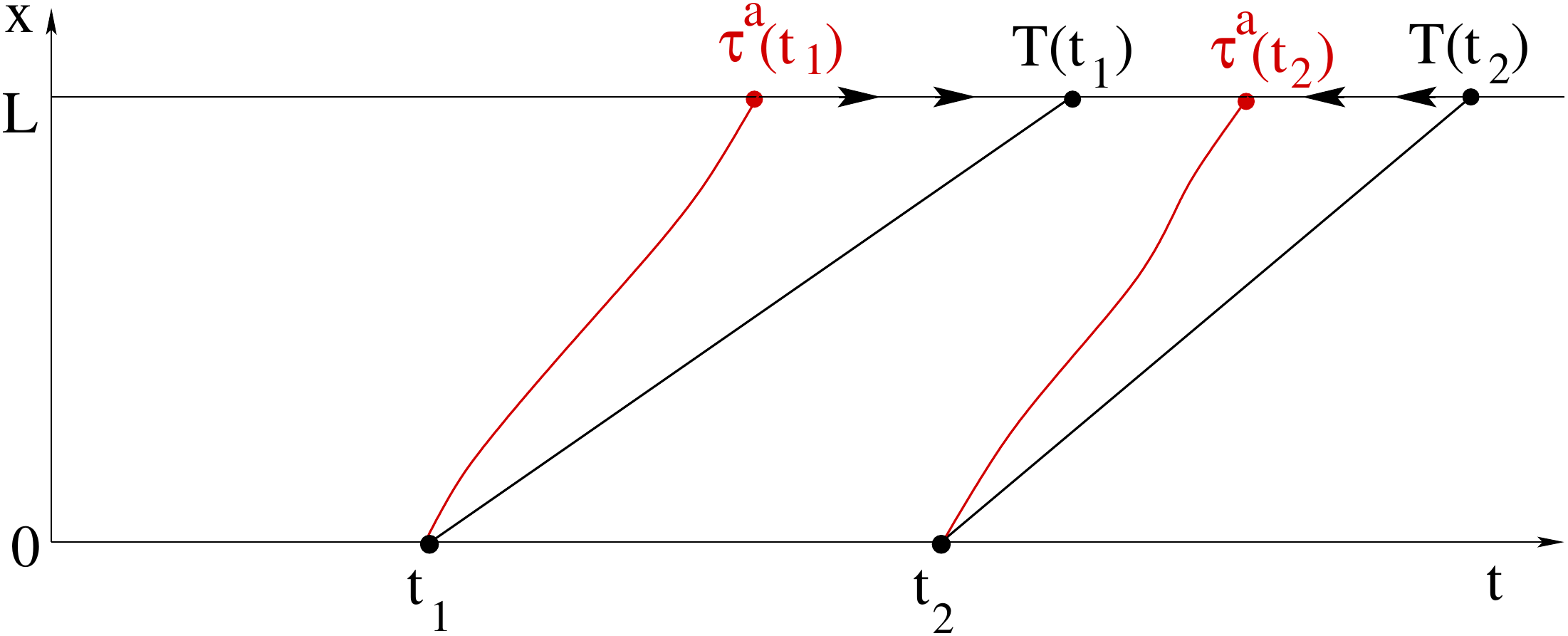}
   \caption{\small Characteristic lines and car trajectories.   Here $\tau^a(t_i)$ is the arrival 
   time of a driver departing at time $t_i$, while $T(t_i)$ is the terminal point of 
   the characteristic
   through $(t_i,0)$.}
   \label{f:tf247}
\end{figure}

The next lemma, analyzing various perturbations to an optimal solution $u$,
provides the key step toward the proof of Theorem~\ref{t:1}.

\begin{lemma}\label{l:3}
Let  $(\bar u_1,\ldots, \bar u_N)= (\bar\theta_1 \bar u,\ldots, \bar\theta_N \bar u)$ be optimal departure rates. 
Assume that $t_1, t_2$
are Lebesgue points for all functions $\bar u, \bar\theta_1,\ldots,\bar\theta_N$, 
and
\bel{a3}\bar u(t_1)\,<\,M,\qquad\qquad \bar u_i(t_2)\,>\,0.\eeq
Call $\tau^a(t_1), \tau^a(t_2)$ the arrival times of a driver departing at times
$t_1, t_2$, respectively. Moreover, let $T(t_1), T(t_2)$ be the times where the 
(unique) generalized forward characteristic starting from $t_1, t_2$ reaches 
the point $L$.  Then 
\bel{nc1}\bega{l}\ds
\vp(t_1) + \psi_i(\tau^a(t_1)) +\sum_{j=1}^N\int_{\tau^a(t_1)} ^{T(t_1)} 
\psi_j'(s) \theta_j(s,L)\, ds \\[4mm]
\qquad\ds 
\geq ~\vp(t_2) + \psi_i(\tau^a(t_2)) +\sum_{j=1}^N\int_{\tau^a(t_2)} ^{T(t_2)}\psi_j'(s) \theta_j(s,L)\, ds \,.\enda
\eeq
\end{lemma}
\v
\begin{remark} {\rm The left hand side of (\ref{nc1}) can be interpreted as the 
cost for inserting an additional $i$-driver, departing at time $t_1$.
In this case, an additional driver arrives at time $T(t_1)$, but this is not the same one!
Indeed, the new driver arrives at time $\tau^a(t_1)$. However, the presence of 
this additional car slows down  all the other cars whose arrival time is  
$T\in [\tau^a(t_1), \, T(t_1)]$. The delay in the arrival time of all these cars 
causes a further  increase in the total cost, accounted by the integral term 
on the left hand side of  (\ref{nc1}).   Similarly, the right hand side
is the amount which can be saved by removing an $i$-driver
departing at time $t_2$.}
\end{remark}

{\bf Proof of Lemma~\ref{l:3}.}   
{\bf 1.}
Since $t_1, t_2$ are Lebesgue points of $\bar u$, they cannot be the center of a rarefaction wave.   Hence  there exist unique points 
$T_1=T(t_1)$ and $ T_2=T(t_2)$ such that 
\bel{t1u}t_i~=~\argmin_\tau 
~\left\{ L g^*\Big( {T_i-\tau\over L}\Big) + \ov U(\tau)\right\},\qquad\qquad i=1,2.\eeq
Assuming that (\ref{nc1}) fails, we shall derive a contradiction.
Indeed, we will construct a new initial data $\bar u_i^\dag$ 
which is slightly smaller than $\bar u_i$ in a neighborhood
of $t_1$ and slightly larger than $\bar u_i$ in a neighborhood of $t_2$,
 yielding a lower total cost.
 Various cases can arise, depending on the relative position of 
 $\tau^a(t_i)$ and $T(t_i)$.  To fix the ideas, in the following we assume that
 \bel{tTi}
 \tau^a(t_1)~<~T(t_1)~<~\tau^a(t_2)~<~T(t_2),\eeq
as shown in Fig.~\ref{f:tf247}.   The other cases are handled in a similar way.

\begin{figure}[htbp]
   \centering
 \includegraphics[width=0.9\textwidth]{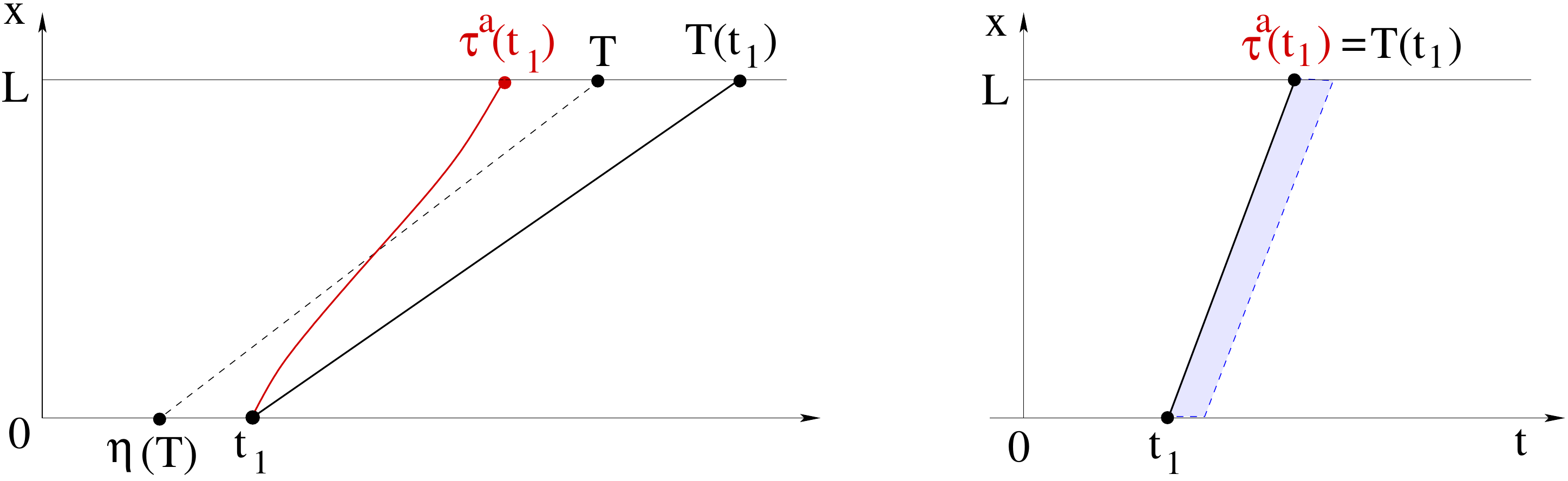}
   \caption{\small Left: if $\tau^a(t_1)<T(t_1)$, the function $u(\cdot,L)$ must be strictly
   positive on the entire interval $[\tau^a(t_1),T(t_1)]$.  
   Right: if $\tau^a(t_1)=T(t_1)$, then the characteristic
   and the car trajectory starting at $(t_1,0)$ coincide, and 
   $u=0$ along this line.  A small perturbation in the initial data, supported on $[t_1, t_1+\delta]$, will modify the solution only in a small neighborhood of this characteristic.
   }
   \label{f:tf249}
\end{figure}

We observe that the above strict inequalities imply that 
$u(\cdot,L)$ is strictly positive on the intervals $[\tau^a(t_1), T(t_1)]$
and $[\tau^a(t_2),T(t_2)]$. Indeed, the car speed is always $\leq v(0)= {1\over g'(0)}$.
As shown in Fig.~\ref{f:tf249}, if 
$${\tau^a(t_1)-t_1\over L}~=~{1\over v(0)}~=~g'(0),$$
then the car speed would be identically equal to the maximum speed $v(0)$.
In this case the car trajectory coincides with a characteristic line, and 
hence $T(t_1)=\tau^a(t_1)$,
against the assumption (\ref{tTi}).   Therefore, we must have
$${\tau^a(t_1)-t_1\over L}~>~g'(0)\,,\qquad\qquad 
\tau^a(t_1) - {L\over v(0)} ~<~t_1\,.$$
Since characteristics do not cross each other, for every
$T\in \,]\tau^a(t_2),T(t_2)[\,$ the initial point of a characteristic through $(T,L)$ must
satisfy 
$$\eta(T)~=~T- L \cdot g'(u(T,L))~\leq~t_1\,.$$
Hence 
\bel{g'u}g'(u(T,L))~\geq~ {T-t_1\over L}~\geq~ {\tau^a(t_1)-t_1\over L}>~g'(0).\eeq
Since $g'$ is an increasing function, this yields a lower bound on $u(T,L)$.

\v
{\bf 2.} %Let  $\ve>0$ be given, and let $\delta, \delta'>0$ be as in 
Consider a perturbed  set of initial data of the form
$(\bar u_1,\ldots, \bar u_i^\dagger, \ldots,\bar u_N)$, where only the 
component $\bar u_i$ is modified.   The new departure rate for 
drivers of the $i$-th group is chosen so that
\bel{ub1}\left\{
\bega{cl} \bar u_i^\dagger(t)~=~\bar u_i(t)\qquad\qquad &\hbox{if}~~t
\notin [t_1\,,~t_1+\delta']\cup[t_2\,,~t_2+\delta']\,,\\[3mm]
\bar u_i^\dagger(t)~\geq~\bar u_i(t)\qquad\qquad &\hbox{if}~~t
\in [t_1\,,~t_1+\delta']\,,\\[3mm]
0\,\leq \,\bar u_i^\dagger(t)~\,\leq\,\bar u_i(t)\qquad\qquad &\hbox{if}~~t
\in [t_2\,,~t_2+\delta']\,,\\[3mm]
\enda\right.\eeq
\bel{ub2}\int_{t_1}^{t_1+\delta'} \Big[\bar u^\dagger_i(s)-\bar u_i(s)
\Big]\, ds
~=~\int_{t_2}^{t_2+\delta'}  \Big[\bar u_i(s)-\bar u^\dagger_i(s)
\Big]\, ds
~=~\delta~>~0\,.\eeq
Given $\ve>0$, according to Lemma~\ref{l:5}, we can choose $\delta,\delta'>0$
small enough so that the perturbation in the initial datum
$$u(t,0)~=~\bar u(t)~=~\sum_{i=1}^N \bar u_i(t)$$
affects the values of $u(\cdot,L)$
only in a small neighborhood of the points $T(t_1), T(t_2)$, namely
\bel{tcm}
 u^\dagger (t,L)~=~ u(t,L) \qquad\qquad\forall t\notin  [T(t_1)-\ve\,,~T(t_1)+\ve]
\cup [T(t_2)-\ve\,,~T(t_2)+\ve].\eeq
\v
{\bf 3.}  We now consider a sequence of perturbations of the form (\ref{ub1})-(\ref{ub2}),
with $\ve_n, \delta_n,\delta'_n\to 0$.
Calling $\tau^a_n(t)$ the corresponding arrival times,
 we claim that the following holds.
 \begi
 \item[{\bf (C)}]
{\it 
Let  
$t$ be a 
Lebesgue point for $\bar u$, with $\bar u(t)>0$, and let $\tau^a(t)$ be a Lebesgue
point for $u(\cdot,L)$. Then $u(\tau^a(t),L)>0$ and  the following  implications hold.
\bel{li1}\tau^a(t_1)~<~\tau^a(t)~<~T_1\qquad \implies\qquad\lim_{n\to\infty}
{\tau^a_n(t)-\tau^a(t)\over\delta_n}~=~{1\over u(\tau^a(t),L)}\,,\eeq
\bel{li2}\tau^a(t_2)~<~\tau^a(t)~<~T_2\qquad \implies\qquad\lim_{n\to\infty}
{\tau^a_n(t)-\tau^a(t)\over\delta_n}~=~-{1\over u(\tau^a(t),L)}\,,\eeq
\bel{li3}\tau^a(t)~\notin~[\tau^a(t_1),T_1]\cup [\tau^a(t_2),T_2]\qquad 
\implies\qquad\lim_{n\to\infty}
{\tau^a_n(t)-\tau^a(t)\over\delta_n}~=~0.\eeq
}
\endi
In first approximation, the above limits show that:
\begi
\item For those drivers who were reaching destination 
at a time $T\in \,]\tau^a(t_1), T_1[\,$,
the arrival time is delayed by $\delta_n/u(\tau^a(t),L)$.
\item For those drivers who were reaching destination 
at a time $T\in \,]\tau^a(t_2), T_2[\,$,
the arrival time is anticipated by $\delta_n/u(\tau^a(t),L)$.
\item For all other drivers, the arrival time does not change.
\endi
To prove the above claim we first observe that, if $u(\tau^a(t),L)=0$, then 
$u(t', x')=0$ along the backward characteristic
$$\{(t',x')\,;~~x'~=~L- (\tau^a(t)-t') v(0)\}.$$
but in this case, this characteristic coincides with a car trajectory.  Hence
$\bar u(t)=u(0, t)=0$ as well, contradicting our first assumption.

To prove (\ref{li1}), assume $\tau^a(t_1)<\tau^a(t)<T_1$. 
Then, for all $n$ large enough, the arrival time $\tau^a_n(t)$ is uniquely determined by the
identity
\bel{IU+}
U(\tau_n^a(t),L)~=~U(\tau^a(t),L) + \delta_n\,.\eeq
Observing that the partial derivative is
$${\partial\over\partial\tau} U(\tau,L)\bigg|_{\tau=\tau^a(t)} ~=~u(\tau^a(t),L),$$
from (\ref{IU+}) one obtains (\ref{li1}).   Notice that here the denominator is uniformly positive, as a consequence of (\ref{g'u}).

The proof of (\ref{li2}) is entirely similar,
replacing (\ref{IU+}) with  the identity
\bel{IU-}
U(\tau_n^a(t),L)~=~U(\tau^a(t),L) - \delta_n\,.\eeq
Finally, if the condition on the left hand side of (\ref{li3}) holds, then for
all $n\geq 1$ sufficiently large one has $\tau_n^a(t)=\tau^a(t)$, and the implication is trivial.
\v
{\bf 4.} By the properties (\ref{li1}) it follows
$$\bega{l}\ds
\lim_{n\to\infty}~{1\over\delta_n}  \sum_{j=1}^N
 \int_{\tau^a(t)\in \,]\tau^a(t_1), T_1[}[\psi_j(\tau^a_n(t))-\psi_j(\tau^a(t))]\, \bar u_j(t)\, dt
 \\[4mm]
 \qquad =~\ds  \sum_{j=1}^N
\int_{\tau^a(t)\in \,]\tau^a(t_1), T_1[}\psi'_j(\tau^a(t))\cdot \lim_{n\to\infty}
{\tau_n^a(t)-\tau_a(t)\over \delta_n}\, \bar u_j(t)\, dt\\[4mm]
 \qquad =~\ds \sum_{j=1}^N
\int_{\tau^a(t_1)}^{ T_1}\psi'_j(\tau) \cdot {1\over u(\tau,L)} \,   u_j(\tau,L)\, d\tau\\[4mm]
\qquad =~\ds \sum_{j=1}^N
\int_{\tau^a(t_1)}^{ T_1}\psi'_j(\tau)  \,   \theta_j(\tau,L)\, d\tau\,.
\enda $$
An entirely similar computation can be performed on the interval $]\tau^a(t_2), T_2[\,$.
Combining these estimates, we thus conclude
\bel{DJ}\bega{l}\ds
\lim_{n\to\infty}~ {J(\bar u_1,\ldots,\bar u_i^\dagger,\ldots, \bar u_N) - J(\bar u_1,\ldots,\bar u_i,\ldots, \bar u_N) \over \delta_n}\\[4mm]
\qquad =~\ds 
\vp(t_1) + \psi_i(\tau^a(t_1)) +\sum_{j=1}^N\int_{\tau^a(t_1)} ^{T(t_1)} 
\psi_j'(s) \theta_j(s,L)\, ds\\[4mm]
\qquad\ds -
\left[\vp(t_2) + \psi_i(\tau^a(t_2)) +\sum_{j=1}^N\int_{\tau^a(t_2)} ^{T(t_2)}\psi_j'(s) \theta_j(s,L)\, ds\right].
\enda\eeq
If the inequality (\ref{nc1}) does not hold, then the right hand side of (\ref{DJ})
is negative.  This yields a contradiction with the optimality of the departure rates 
$(\bar u_1,\ldots,\bar u_i,\ldots, \bar u_N)$.
\v
{\bf 5.} The above analysis proves the lemma in the case where (\ref{tTi}) holds.
On the other hand, if $T(t_1)= \tau^a(t_1)$, then a small perturbation of the 
departure rate on the interval  $[t_1, t_1+\delta_n]$ will modify the arrival rate
only in a small neighborhood of $\tau^a(t_1)$ (see Fig.~\ref{f:tf249}, right).   In this case, one directly proves that 
the limit (\ref{DJ}) remains valid, since the  integral over
the interval $[\tau^a(t_1), \,T(t_1)]$ trivially vanishes.
\endproof
\v
\subsection{Proof of Theorem~\ref{t:1}.}

Let $u_i(t,0)=\bar\theta_i(t)\bar u(t)$ be  optimal departure rates, 
and let $u$, $\theta_i$   be the corresponding solutions to 
\bel{cli} 
u_x + g(u)_t~=~0,\qquad\qquad \theta_{i,t} + v(g(u)) \theta_{i,x}
~=~0,
\qquad i=1,\ldots, N.
\eeq
%We observe that necessary conditions for optimality will be derived by means of two types of
%perturbations:
%\begi
%\item[(i)] Changing the order in which drivers of different group 
%depart,  but without modifying the overall departure rate $\bar u$.
%This type of perturbations do not affect the solution $u=u(t,x)$.

%\item[(ii)] Changing the departure rate $\bar u$, but without modifying the 
%order in which different drivers depart.    This type of perturbations do not affect the 
%functions $\Theta_i(\cdot)$ in (\ref{TIS}).
%\endi
%\v
The proof will be worked out in several steps.
\v
{\bf 1.} We begin by showing that an optimal solution no shock can occur
in the interior of the domain, i.e.~for $0\leq x<L$.

Indeed, assume on the contrary that a shock is present, and let $(T,L)$ be the terminal
position of this shock.  According to Remark~\ref{r:shock}, we can 
change the initial datum so that the new solution $u^\dagger$ 
contains a centered compression 
wave focusing at $(T,L)$.   
More precisely,  define $\ov U^\dagger$, $\bar u^\dagger$ 
as in Remark~\ref{r:shock}.
Assuming that the functions $\Theta_i(s)$ are defined by
\bel{ovu}
\ov U_i(t)~\doteq~\Theta_i(\ov U(t))\cdot \ov U(t),\eeq
 define the components
$ (\ov U^\dagger_1,\ldots, \ov U^\dagger_N)$ by setting
\bel{ovui}
\ov U_i^\dagger(t)~\doteq~\Theta_i(\ov{U}^{\dagger}(t))\cdot \ov{U}^{\dagger}(t).\eeq
Notice that these definition imply
$$U_i^\dagger(t,L)~=~U_i(t,L)\qquad\qquad\forall t\in\R,$$
hence the arrival costs remain the same. On the other hand, we have
\bel{OV3}
\ov U_i^\dagger(t)~\leq~\ov U_i(t),\eeq
for all $i,t$.
Observing that there exists some $i$ and some 
$t\in   [\eta^-(T), \,\eta^+(T)]$ where (\ref{OV3}) is satisfied as a strict inequality,
we claim that the total departure
cost  for the perturbed solution is strictly smaller.

To see this, introduce a variable $\xi\in [0, G_i]$ labeling drivers of the $i$-th group.
Define the departure times 
\bel{tix}
t_i(\xi)~\doteq~\inf\{t\,;~~\ov U_i(t)>\xi\},\qquad\quad
t^\dagger_i(\xi)~\doteq~\inf\{t\,;~~\ov U^\dagger_i(t)>\xi\}.\eeq
By the previous definitions it follows
$$t_i(\xi)~\leq~t_i^\dagger(\xi),$$
with strict inequality holding at least for some index $i$ and some values of $\xi\in [0, G_i]$.
We now compute
$$\sum_i \int \vp(t) d \ov U_i^\dagger (t)~=~\sum_i \int_0^{G_i}
\vp(t_i^\dagger(\xi)) \, d\xi~<~\sum_i \int_0^{G_i}
\vp(t_i(\xi)) \, d\xi~=~\sum_i \int \vp(t) d \ov U_i (t),
$$
proving our claim.
\v
{\bf 2.} Next, we claim that an optimal departure rate satisfies
\bel{BUM}\bar u(t)~<~M\qquad\quad\forall t\in\R.\eeq
Indeed, since the characteristic speed 
satisfies $g'(u)\to +\infty$ as $u\to M$, 
if $\bar u(\tau)=M$ at some point $\tau$ then the solution $u(\cdot,x)$
would immediately contain a shock, for every $x>0$.   By the previous step,
this contradicts the optimality assumption.
\v
{\bf 3.} 
According to Lemma~\ref{l:3}, by (\ref{BUM}), the quantity
\bel{DJ2}\Delta J(i,t)~=~\vp(t) + \psi_i(\tau^a(t)) +\sum_{j=1}^N\int_{\tau^a(t)}^{T(t)} 
\psi_j'(s)\, \theta_j(s,L)\, ds
\eeq
is equal to some constant $C_i$ for all $t\in {\rm Supp}(\bar u_i)$, and is greater or equal to
$C_i$ for all $t\in\R$.  In other words, for each $i=1,\ldots,N$ we have
\bel{J2}\vp(t) + \psi_i(\tau^a(t)) - C_i +\sum_{j=1}^N\int_{\tau^a(t)}^{T(t)} 
\psi_j'(s)\, \theta_j(s,L)\, ds~=~0\qquad\forall t\in \hbox{Supp}(\bar u_i)\,,\eeq
\bel{J3}\vp(t) + \psi_i(\tau^a(t)) - C_i +\sum_{j=1}^N\int_{\tau^a(t)}^{T(t)} 
\psi_j'(s)\, \theta_j(s,L)\, ds~\geq ~0\qquad\forall t\in \R.\eeq
This implies 
\bel{J4} \vp(t) + \psi_i(\tau^a(t)) - C_i ~=~\vp(t) +\min_j \Big(\psi_j(\tau^a(t)) - C_j\Big)~=~
-\sum_{j=1}^N\int_{\tau^a(t)}^{T(t)} 
\psi_j'(s)\, \theta_j(s,L)\, ds\eeq
for all $ t\in \hbox{Supp}(\bar u_i)$.
We now observe that, for a.e.~$s\in [\tau^a(t), T(t)]$ and $j\in \{1,\ldots,N\}$, one has the implication
\bel{thp}\theta_j(s,L)\,>\,0\qquad\implies\qquad \psi_j(s)-C_j~=~\min_k (\psi_k(s)-C_k).\eeq
Moreover, for every $j,k$ and a.e.~$s$ in the set 
$$\bigl\{s\in\R\,;~~\psi_j(s)-C_j = \psi_k(s)-C_k\bigr\},$$ one has
$\psi_j'(s) = \psi_k'(s)$. Therefore, defining $\psi$ as in (\ref{pps})
and recalling that  $\sum_j\theta_j=1$, we obtain
\bel{J7} \sum_{j=1}^N\int_{\tau^a(t)}^{T(t)} 
\psi_j'(s)\, \theta_j(s,L)\, ds~=~\int_{\tau^a(t)}^{T(t)} 
\psi'(s)\, ds.\eeq
In turn, this implies
\bel{J5}\bega{l}\ds\vp(t) + \psi_i(\tau^a(t)) - C_i +\sum_{j=1}^N\int_{\tau^a(t)}^{T(t)} 
\psi_j'(s)\, \theta_j(s,L)\, ds\\[4mm]
\qquad =~\ds
\vp(t) + \psi_i(\tau^a(t))- C_i+ \int_{\tau^a(t)}^{T(t)} 
\psi'(s)\, ds\\[4mm]
\qquad =~\ds \vp(t) + \psi_i(\tau^a(t)) - C_i +\psi(T(t))- \psi(\tau^a(t))\\[4mm]
\qquad =~\ds\vp(t) + \psi(T(t))~=~0.\enda\eeq
According to (\ref{J5}), each characteristic where the solution $u$ is positive 
must connect two points
$(t,0)$ with $(T(t), L)$ with $\vp(t)+\psi(T(t))=0$.
This proves part (i) of Theorem~\ref{t:1}.  Finally, part (ii) follows from (\ref{thp}).
\endproof

\section{An algorithm to construct optimal solutions}  
\label{s:5}
\setcounter{equation}{0}

Here
we illustrate how these necessary conditions can be used to 
construct optimal solutions.  For simplicity, we shall assume that the cost functions $\psi_i$ 
are $\C^2$ and satisfy the assumption
\begi\item[{\bf (A3)}] For any $i\not= j$ one has the implication
\bel{gas}
\psi_i'(t)~=~\psi_j'(t)\qquad\implies\qquad \psi''(t)\not=\psi''_j(t).\eeq
\endi
Notice that, by (\ref{gas}),  for any given constants $C_1, \ldots, C_N$, the set of times
$$\{t\in\R\,;~~\psi_i(t)-C_i~=~\psi_j(t)-C_j\qquad\hbox{for some }~i\not= j\}$$
consists only of isolated points, hence it has measure zero.

We remark that the assumption {\bf (A3)} is generically valid in the space of 
twice continuously differentiable functions.  
Indeed, given $\ve>0$ and any $N$-tuple of twice continuously differentiable 
functions $(\Hat \psi_1,\ldots,\Hat \psi_N)$, by a small perturbation
one can construct functions $\psi_1,\ldots,
 \psi_N$ which satisfy {\bf (A3)} together with 
$$\|\Hat \psi_i-\psi_i\|_{\C^2}~<~\ve,\qquad\qquad i=1,\ldots, N.$$

Let now $G_1,\ldots, G_N$ be the sizes of the $N$ groups of drivers.
In order to construct a globally optimal  
family of departure rates $u_1(\cdot),\ldots, u_N(\cdot)$, we introduce 
the following algorithm.
\begi
\item[(i)] Start by guessing $N$ constants $C_1,\ldots, C_N$, and define the 
 cost function $\psi$ as in (\ref{pps}).
\item[(ii)] Let  $u=u(t,x)$ be the solution of (\ref{CP}) constructed according to
 (\ref{ncc})-(\ref{uc}).
\item[(iii)] Define the sets
\bel{Aidef}A_i~=~\Big\{t\in\R\,;~~u(t,L)>0,\quad \psi_i(t)-C_i \,= \,
\min_k (\psi_k(t)-C_k)\Big\}.\eeq 
Notice that, by {\bf (A3)}, for a.e.~$t\in\R$ the minimum 
in (\ref{Aidef}) is attained by a unique index $k\in \{1,\ldots,N\}$.
\item[(iv)] Consider the map 
$\Lambda:\R^N\mapsto\R^N$, ~   $\Lambda (C_1,\ldots, C_N)
= (\kappa_1,\ldots, \kappa_N)$, where $\kappa_i$ is the total number  
of drivers of the $i$-th group, defined by
\bel{tfi}
\kappa_i~\doteq~\int_{A_i} u(t,L)\, dt\,.\eeq
Determine values $ \Hat C_1,\ldots,  \Hat C_N$ such that 
\bel{LCG}\Lambda ( \Hat C_1,\ldots,  \Hat C_N)
~=~ (G_1,\ldots, G_N).\eeq
\item[(v)]
Set $$ \psi (t)~= ~\min_i ~(\psi_i(t) - \Hat C_i)\,.$$  
And let $u= u (t,x)$ be the solution of (\ref{CP}) whose characteristics satisfy (\ref{ncc}).
\item[(vii)] Finally, for $T\in A_i$, call
$\eta(T)$ the departure time of the driver that arrives at time $T$.
This is obtained by solving the ODE
\bel{ct}\dot x(t)~=~v(t, x(t))~=~{u(t,x)\over g(u(t,x))}\,,\eeq
with terminal condition
\bel{tcT}x(T)=L\,.\eeq
The solution $t\mapsto x(t,T)$ of (\ref{ct})-(\ref{tcT}) yields the trajectory of a car
arriving at the end of the road time $T$.  Its departure time
$\eta(T)$ is defined by the equality $x(\eta(T))=0$.

We now consider the sets of departure times
$$A_i^* ~\doteq~\{\eta(t)\,;~~t\in A_i\}.$$
The departure distribution
$$\bar u_i(t)~\dot=~\left\{ \bega{cl} \bar u(t)\qquad &\hbox{if}\quad t\in A_i^*,\\[3mm]
0\qquad &\hbox{if}\quad t\notin A_i^*,\enda\right.
$$
then satisfies all the necessary conditions for optimality.
\endi
We remark that these conditions are only necessary, not sufficient for optimality.
Since an optimal solution exists, and can obtained by the above method, 
the previous analysis implies that, if the $N$-tuple $(\Hat C_1,\ldots, \Hat C_N)$
which satisfies (\ref{LCG}) is unique, then this must yield the optimal solution.

\begin{figure}[htbp]
 \centering
\includegraphics[width=1.0\textwidth]{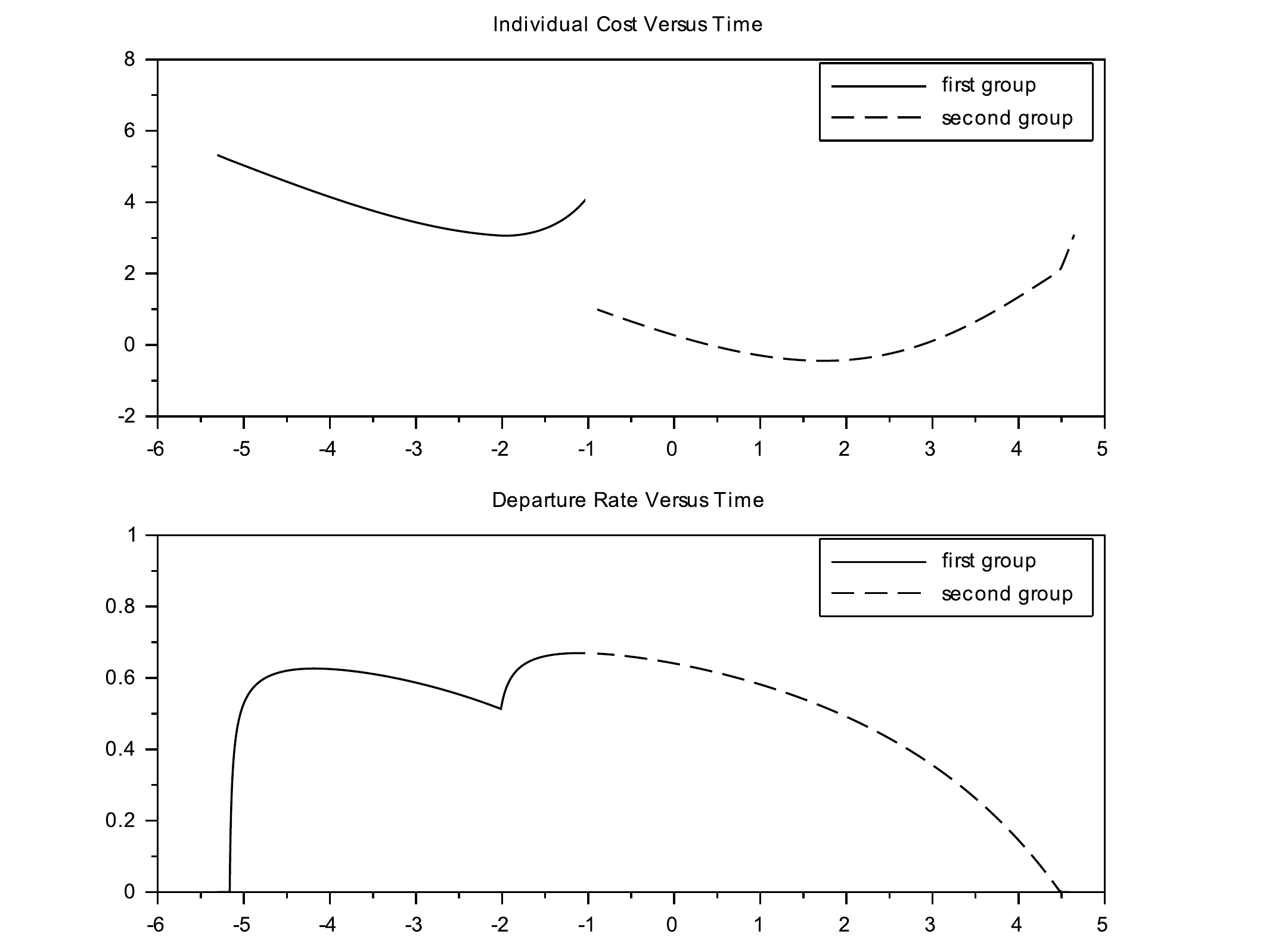}
  \caption{\small The globally optimal solution for the problem described in Example~4. 
  Top: the cost
  incurred by a driver of the first and of the second group, departing at time $t$.
  Bottom: the rate of departure of drivers of the first and of the second group, 
  as a function of time.  
  }
 \label{f:optraf2}
\end{figure}

{\bf Example 4.}
We seek a globally optimal departure rate for 
two groups of drivers,  
with sizes 
$G_1= G_2=2.51$
 on a road with length $L=10$.
The conservation law governing traffic density is
\bel{claw8}\rho_t +[\rho v(\rho)]_x~=~0, \qquad v(\rho) \, =\, 2-\rho\eeq
The departure and arrival costs for drivers of the two groups are
\bel{cost2}\vp(t)~=~-t\qquad \psi_1(t)~=~e^{t-4}\qquad
\psi_2(t)~=~e^{t-7.6}.\eeq
The optimal solution, found by the algorithm described above, 
is shown in Fig.~\ref{f:optraf2}.   The marginal costs for adding one more driver of the first group or of the second group
are
found to be $\Hat C_1=5.18$ and $\Hat C_2 =  2.10$, respectively. 
\v
{\bf Acknowledgment.} This research  was partially supported
by NSF, with grant  DMS-1411786: ``Hyperbolic Conservation Laws and Applications". 
\v

\end{document}